\documentclass{elsarticle}

\usepackage{amssymb,amsthm,amsmath,tikz}
\usepackage{hyperref}

\newtheorem{thm}{Theorem}[section]
\newtheorem{prop}[thm]{Proposition}
\newtheorem{cor}[thm]{Corollary}
\newtheorem{lem}[thm]{Lemma}

{\theoremstyle{definition}
\newtheorem{rem}[thm]{Remark}
\newtheorem{exa}[thm]{Example}
}

\numberwithin{equation}{section}
\numberwithin{figure}{section}

\newcommand{\ble}{\begin{lem}}
\newcommand{\ele}{\end{lem}}
\newcommand{\bth}{\begin{thm}}
\renewcommand{\eth}{\end{thm}}
\newcommand{\bpr}{\begin{prop}}
\newcommand{\epr}{\end{prop}}
\newcommand{\bco}{\begin{cor}}
\newcommand{\eco}{\end{cor}}
\newcommand{\bex}{\begin{exa}}
\newcommand{\eex}{\end{exa}}
\newcommand{\barr}{\begin{array}}
\newcommand{\earr}{\end{array}}
\newcommand{\beq}{\begin{equation}}
\newcommand{\eeq}{\end{equation}}
\newcommand{\bea}{\begin{eqnarray*}}
\newcommand{\eea}{\end{eqnarray*}}
\newcommand{\bprf}{\begin{proof}}
\newcommand{\eprf}{\end{proof}\medskip}
\newcommand{\hqed}{\hfill \qed}

\newcommand{\sbe}{\subseteq}
\newcommand{\spn}[1]{\langle{#1}\rangle}
\newcommand{\ree}[1]{(\ref{#1})}
\newcommand{\ra}{\rightarrow}  
\newcommand{\lra}{\longrightarrow}
\newcommand{\si}{\sigma}
\newcommand{\La}{\Lambda}
\newcommand{\bbP}{{\mathbb P}}
\newcommand{\cI}{{\cal I}}
\newcommand{\fS}{{\mathfrak S}}

\newcommand{\poset}{P}
\newcommand{\pzero}{P_0}
\newcommand{\pstar}{\poset^*}
\newcommand{\len}[1]{|#1|}
\newcommand{\mupstar}{\mu}
\newcommand{\mulambdastar}{\mu}
\newcommand{\mup}{\mu_{0}}
\newcommand{\emb}{\eta}
\newcommand{\covers}{\lra}
\newcommand{\coversl}[1]{\stackrel{#1}{\covers}}
\newcommand{\coversinp}{\lra_0}
\newcommand{\chainlabel}[1]{L(#1)}
\newcommand{\lexleq}{\leq_{\mathrm{lex}}}

\newcommand{\lexgreater}{>_{\mathrm{lex}}}
\newcommand{\lexless}{<_{\mathrm{lex}}}
\newcommand{\zeros}{\mathrm{zero}}
\newcommand{\czero}{C_0}
\newcommand{\ptomie}{\Lambda_s}
\newcommand{\muembpstar}[2]{\mupstar_{\mathrm{emb}}(#1,#2)}
\newcommand{\jint}{\mathcal{J}}
\newcommand{\arbchain}{\mathfrak{C}}
\newcommand{\defect}{\mathrm{def}}
\DeclareMathOperator{\rk}{rk}

\newcommand{\wdots}{\cdots}
\definecolor{boxes}{RGB}{0,0,255}

\journal{Advances in Mathematics}

\begin{document}

\begin{frontmatter}

\title{The M\"obius function of generalized subword order}

\author{Peter R. W. McNamara\fnref{peter}}
\ead{peter.mcnamara@bucknell.edu}
\ead[url]{http://www.facstaff.bucknell.edu/pm040}

\author{Bruce E. Sagan\fnref{bruce}}
\ead{sagan@math.msu.edu}
\ead[url]{http://www.math.msu.edu/~sagan}

\address[peter]{Department of Mathematics, Bucknell University, Lewisburg, PA 17837, USA}

\address[bruce]{Department of Mathematics, Michigan State University, East Lansing, MI 48824-1027, USA}

\begin{abstract}
Let $P$ be a poset and let $P^*$ be the set of all finite length words over $P$.  Generalized subword order is the partial order on $P^*$ obtained by letting $u\le w$ if and only if there is a subword $u'$ of $w$ having the same length as $u$ such that each element of $u$ is less than or equal to the corresponding element of $u'$ in the partial order on $P$.  Classical subword order arises when $P$ is an antichain, while letting $P$ be a chain gives an order on compositions.  For any finite poset $P$, we give a simple formula for the M\"obius function of $P^*$ in terms of the M\"obius function of $P$.  This permits us to rederive in a easy and uniform manner previous results of Bj\"orner, Sagan and Vatter, and Tomie.  We are also able to determine the homotopy type of all intervals in $P^*$ for any finite $P$ of rank at most 1.
\end{abstract}


\end{frontmatter}

\section{Introduction}

Let $A$ (the {\it alphabet\/}) be any set  and let $A^*$ be the {\it Kleene closure\/} of all finite length words over $A$, so
$$
A^*=\{w=w(1)w(2)\wdots w(n) : \mbox{$w(i)\in A$ for all $i$, and $n\ge0$}\}.
$$
We denote the \emph{length} or \emph{cardinality} of $w$ by $\len{w}$.
A {\it subword\/} of $w\in A^*$ is a word $u=w(i_1)w(i_2)\wdots w(i_k)$ where $i_1<i_2<\dots<i_k$.  (Note that the elements chosen from $w$ need not be consecutive.)  {\it Subword order\/} on $A^*$ is defined by letting $u\le w$ if and only if $u$ is a subword of $w$.  Bj\"orner~\cite{bjo:mfs} was the first person to determine the M\"obius function of subword order.

Now consider the symmetric group $\fS_n$ of all permutations of $\{1,2,\ldots,n\}$.  If $\si=\si(1)\si(2)\wdots\si(n)\in\fS_n$ and $\pi=\pi(1)\pi(2)\wdots\pi(k)\in\fS_k$ then {\it $\si$ contains a copy of $\pi$ as a pattern\/} if there is a subword $\si(i_1)\si(i_2)\wdots\si(i_k)$ such that 
$$
\pi(r)<\pi(s) \iff \si(i_r)<\si(i_s)
$$
for all $1\le r<s\le k$.  The {\it pattern order\/} on $\fS=\uplus_{n\ge0}\fS_n$ is obtained by letting $\pi\le\si$ if and only if $\si$ contains a copy of $\pi$.  For example, $2143 \leq 321465$ because of the subwords $3265$, $3165$ or $2165$.  Wilf \cite{wil:pp} posed the problem of determining the M\"obius function of pattern order.  The first result along these lines was obtained by Sagan and Vatter~\cite{sv:mfc} and this will be discussed in more detail below.  Later work has been done by Steingr\'{\i}msson and Tenner \cite{st:mfp} and by Burstein, Jel{\'{\i}}nek, Jel{\'{\i}}nkov\'a and Steingr\'{\i}msson \cite{bjjs:mfs}.  
It remains an open problem to fully answer Wilf's question.

When trying to prove results about pattern containment, it is often instructive to consider the case of {\it layered\/} permutations, which are those of the form
$$
\pi = a, a-1,\dots,1,a+b,a+b-1,\dots,a+1,a+b+c,a+b+c-1,\dots
$$
for some positive integers $a,b,c,\dots$.  Note that a layered permutation is completely specified by the composition $(a,b,c,\dots)$ of {\it layer lengths\/}, and that pattern order on layered permutations is isomorphic to the following order on compositions: for compositions $a = (a_1, a_2, \dots, a_r)$ and $b=(b_1, b_2, \dots, b_s)$, we say that $a \leq b$ if there exists a subsequence $(b_{i_1}, b_{i_2}, \dots, b_{i_r})$ of $b$ such that $a_j \leq b_{i_j}$ for $1 \leq j \leq r$.  Our example $2143 \leq 321465$ above for layered permutations corresponds to $22 \leq 312$ for compositions.

\begin{figure}
\begin{center}
\begin{tikzpicture}
\tikzstyle{elt}=[rectangle, draw=boxes]
\matrix{
&\node(3)[elt]{$3$};&\\[20pt]
\node(1)[elt]{$1$};&[10pt]&[10pt]\node(2)[elt]{$2$};\\
};
\draw (3)--(1);
\draw (3)--(2);
\end{tikzpicture}
\caption{The poset $\Lambda$}
\label{Lambda:fig}
\end{center}
\end{figure}
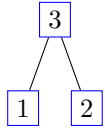

Sagan and Vatter~\cite{sv:mfc} generalized both subword order and pattern order on layered permutations as follows.  
Letting $P$ be any poset, it is natural to let $P^*$ denote the Kleene closure of the alphabet consisting of the elements of $P$.  Define
{\it generalized subword order\/} on $P^*$ by letting $u\le w$ if and only if there is a subword $w(i_1)w(i_2)\wdots w(i_k)$ of $w$  of the same length as $u$ such that 
\beq
\label{gso:def}
u(j) \leq_P w(i_j) \mbox{\ \ for $1\leq j \leq k.$}
\eeq
Note that if $P$ is an antichain, then generalized subword order on $P^*$ is the same as ordinary subword order since one can only have $a\le_P b$ if $a=b$.  At the other extreme, if $P$ is the chain $\bbP$ of positive integers, then, as remarked in the previous paragraph, generalized subword order on the set $\bbP^*$ of compositions is isomorphic to pattern order on layered permutations.   Sagan and Vatter determined the M\"obius function of $P^*$ for any rooted forest $P$, i.e., each component of the Hasse diagram of $P$ is a tree with a unique minimal element.  Note that this covers both the antichain and chain cases.  They also considered the smallest $P$ which is not a rooted forest, namely the poset $\La$ given in Figure~\ref{Lambda:fig}, and conjectured that the M\"obius values for certain intervals in $\La^*$ were given by coefficients of Chebyshev polynomials of the first kind.  This conjecture was later proved and the result generalized by Tomie~\cite{tom:gcp} using ad hoc methods.  Earlier appearances of generalized subword order in the context of well-quasi-orderings are surveyed in \cite{kru:twq}.

Our main result is a simple formula for the M\"obius function of $P^*$ for any finite poset $P$, as given in Theorem~\ref{main:thm}.  To state the theorem, we need to introduce some key notation and terminology.  We let $P_0$ denote the poset $P$ with a bottom element 0 adjoined, and let $\leq_0$ and $\mup$ denote the order relation and M\"obius function of $P_0$, 
 respectively.   An {\it expansion\/} of $u\in P^*$ is $\emb\in (P_0)^*$ such that the restriction of $\emb$ to its nonzero elements is $u$.  For example, $0110302$ is an expansion of $u=1132$.  An {\it embedding of $u$ in $w$\/} is an expansion $\emb$ of $u$ having length $\len{w}$ such that $\emb(j)\leq_0 w(j)$ for all $j$.  Continuing our example and using the poset in Figure~\ref{Lambda:fig}, we see that the given expansion can be considered as an embedding of $u$ in $w=2132333$.  It should be clear from the definitions that there is an embedding of $u$ in $w$ if and only if $u\le w$ in $P^*$.  Since the M\"obius function of $P^*$ is our principal object of interest, we abbreviate $\mu_{P^*}$ by $\mu$.  With these fundamentals in place, we can now state our main result.

\begin{figure}[htb] 
\begin{center}
\hspace*{-5mm}
\begin{tikzpicture}
\tikzstyle{elt}=[rectangle, draw=boxes]
\matrix{
&[15pt]&[15pt]&[15pt]&[15pt]&\node(333)[elt]{$333$};&[15pt]&[15pt]&[15pt]&[15pt]&[15pt]\\[30pt]
\node(133)[elt]{$133$};&&\node(233)[elt]{$233$};&&\node(313)[elt]{$313$};&&
     \node(323)[elt]{$323$};&&\node(331)[elt]{$331$};&&\node(332)[elt]{$332$};\\[80pt]
\node(33)[elt]{$33$};&
\node(113)[elt]{$113$};&
\node(123)[elt]{$123$};&
\node(131)[elt]{$131$};&
\node(132)[elt]{$132$};&&
\node(213)[elt]{$213$};&
\node(231)[elt]{$231$};&
\node(311)[elt]{$311$};&
\node(312)[elt]{$312$};&
\node(321)[elt]{$321$};\\[80pt]
\node(13)[elt]{$13$};&&\node(31)[elt]{$31$};&&\node(111)[elt]{$111$};&&
     \node(112)[elt]{$112$};&&\node(121)[elt]{$121$};&&\node(211)[elt]{$211$};\\[30pt]
&&&&&\node(11)[elt]{$11$};&&&&&\\
};
\draw (333)--(133);
\draw (333)--(233);
\draw (333)--(313);
\draw (333)--(323);
\draw (333)--(331);
\draw (333)--(332);
   \draw (133)--(33);
   \draw (133)--(113);
   \draw (133)--(123);
   \draw (133)--(131);
   \draw (133)--(132);
      \draw (233)--(33);
      \draw (233)--(213);
      \draw (233)--(231);
         \draw (313)--(113);
         \draw (313)--(213);
         \draw (313)--(33);
         \draw (313)--(311);
         \draw (313)--(312);
            \draw (323)--(123);
            \draw (323)--(33);
            \draw (323)--(321);
               \draw (331)--(131);
               \draw (331)--(231);
               \draw (331)--(311);
               \draw (331)--(321);
               \draw (331)--(33);
                  \draw (332)--(132);
                  \draw (332)--(312);
                  \draw (332)--(33);
   \draw (33)--(13);
   \draw (33)--(31);
      \draw (113)--(13);
      \draw (113)--(111);
      \draw (113)--(112);
         \draw (123)--(13);
         \draw (123)--(121);
            \draw (131)--(31);
            \draw (131)--(111);
            \draw (131)--(121);
            \draw (131)--(13);
               \draw (132)--(112);
               \draw (132)--(13);
                  \draw (213)--(13);
                  \draw (213)--(211);
                     \draw (231)--(31);
                     \draw (231)--(211);
                        \draw (311)--(111);
                        \draw (311)--(211);
                        \draw (311)--(31);
                           \draw (312)--(112);
                           \draw (312)--(31);
                              \draw (321)--(121);
                              \draw (321)--(31);
\draw (11)--(13);
\draw (11)--(31);
\draw (11)--(111);
\draw (11)--(112);
\draw (11)--(121);
\draw (11)--(211);
\end{tikzpicture}
\caption{The interval $[11,333]$ of $\pstar$ in the case where $\poset$ is as shown in Figure~\ref{Lambda:fig}}
\label{Lambda11333:fig}
\end{center}
\end{figure}
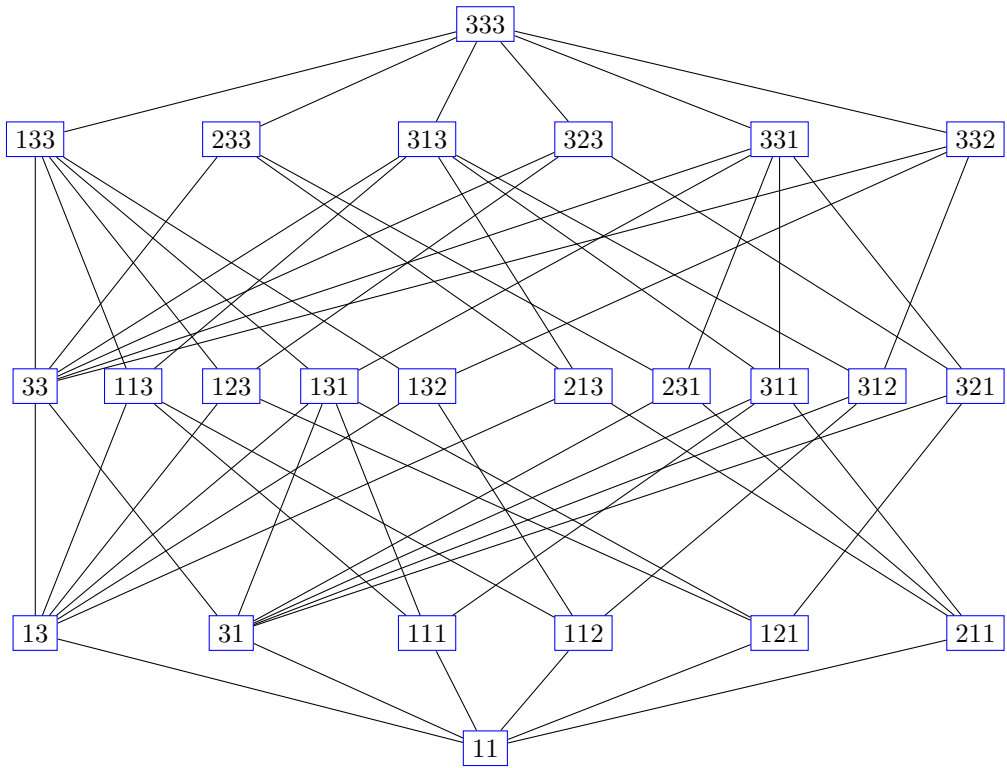

\bth\label{main:thm}
Let $P$ be a poset such that $P_0$ is locally finite.  Let $u$ and $w$ be elements of $\pstar$ with $u \leq w$.  Then
\[
\mupstar(u,w) = \sum_\emb \prod_{1 \leq j \leq \len{w}} \left\{ 
\begin{array}{ll} 
\mup(\emb(j), w(j)) +1 & \mbox{if $\emb(j)=0$ and $w(j-1)=w(j)$}, \\ 
\mup(\emb(j),w(j)) & \mbox{otherwise},
\end{array} \right. 
\] 
where the sum is over all embeddings $\emb$ of $u$ in $w$.  
\eth

If $j=1$, the condition that $w(j-1)=w(j)$ is considered false since $w(j-1)$ does not exist.  The power of Theorem~\ref{main:thm} is that it allows us to determine the M\"obius function in $\pstar$ just by knowing the M\"obius function in $\pzero$; typically, $\pstar$ is a much more complicated poset than $\pzero$, as in Example~\ref{Lambda:exa} below.  We also note that it is natural that a formula for $\mupstar$ involve $\mup$, since if $w$ consists of one letter, then $\mupstar(u,w) = \mup(u,w)$ when $u \leq w$.  

\begin{exa}\label{Lambda:exa}  Let $P=\La$ as shown in Figure~\ref{Lambda:fig}, and consider $\mupstar(11, 333)$. 
\linebreak
 Applying Theorem~\ref{main:thm}, we see that the embedding $\emb=110$ contributes
 \linebreak
 $(-1)(-1)(1+1)=2$ to the sum.  Similarly, $101$ and $011$  contribute 2 and 1, respectively.  Thus $\mupstar(11,333)=5$, which is not at all obvious from Figure~\ref{Lambda11333:fig}.  It is easy to generate intervals whose Hasse diagrams are too large and complicated to be shown clearly here, but whose M\"obius functions are easy to calculate using Theorem~\ref{main:thm}.  One extreme example is that the Hasse diagram of  the interval $[\emptyset, 33333]$ in $\Lambda^*$ has 1904 edges; since the only embedding of the empty word $\emptyset$ in $33333$ is $00000$, Theorem~\ref{main:thm} gives $\mupstar(\emptyset,33333)=(1)(1+1)^4=16$. 
\end{exa}

We prove Theorem~\ref{main:thm} by using Babson and Hersh's method \cite{bh:dmf} for applying Forman's discrete version of Morse theory \cite{for:dmt,for:mtc,for:ugd} to order complexes of posets.  In the next section, we introduce the necessary machinery from \cite{bh:dmf} and the corresponding setup for $P^*$.  Section~\ref{pmr} contains the full proof of Theorem~\ref{main:thm}.
By specializing our result, one can easily derive all the formulas for
M\"obius functions cited above, which we do in Section~\ref{a}.  Specifically, we derive the following results.
\begin{enumerate}
\item Bj\"orner's formula for the M\"obius function of subword order.
\item Sagan and Vatter's result for the M\"obius function of $P^*$ in the case that $P$ is a rooted forest.
\item Their related result for the order on compositions described above, which corresponds to the case when $P = \mathbb{P}$, the positive integers.
\item Tomie's result for the M\"obius function of $\Lambda^*$. The connection to Chebyshev polynomials $T_n(x)$ of the first kind is that $\mu(1^i, 3^j)$ is the coefficient of $x^{j-i}$ in $T_{i+j}(x)$, for all $0 \leq i \leq j$.  
\item Tomie's more general result, which corresponds to letting $P$ consist of an $s$-element antichain with a top element added.  
\end{enumerate}

We can also compute the homotopy type of $P^*$ whenever the rank of $P$, denoted $\rk(P)$, is at most 1; we show that any interval $[u,w]$ in $P^*$ is homotopic to a wedge of $|\mu(u,w)|$ spheres, all of dimension $\rk(w)-\rk(u)-2$.  As a corollary, we get the corresponding result of Bj\"orner~\cite{bjo:mfs} in the antichain case.  The final section contains some concluding remarks about related work.

\section{Generalized subword order and critical chains}
\label{gso}

Here we will only give the minimum amount of detail from the work of Babson and Hersh for the reader to understand our proof of Theorem~\ref{main:thm}.  In particular, we will only talk about the combinatorial application of discrete Morse theory to posets.  Those wishing to find out more about discrete Morse theory itself should consult Forman's excellent primer on the subject \cite{for:ugd}.

\subsection{Discrete Morse theory for posets}

Let $P$ be a locally finite poset, meaning all intervals in $P$ are finite.  If $x,y\in P$, then we write $y\covers x$ if  {\it $y$ covers $x$\/} in that $y>x$ and there is no $z$ with $y>z>x$.  Given an interval $[x,z]$ in $P$, we consider all  (containment) maximal chains $C$ in the interval, so that $C$ must have the form
\beq
\label{C}
C:z=y_0\covers y_1\covers y_2\covers \cdots\covers y_n=x.
\eeq
We wish to totally order the maximal chains in $[x,z]$ in  a way that permits computation of the M\"obius function as in Theorem~\ref{dmt:thm} below.
We say that two maximal chains $C$ and 
$$
C':z=y'_0\covers y'_1\covers y'_2\covers \cdots \covers y'_n=x 
$$
{\it agree to index $k$\/} if 
$$
y_0=y'_0,\ \ y_1=y'_1,\ \ldots\ ,\ y_k=y'_k.
$$
These two chains {\it diverge from index $k$\/} or {\it diverge from the element $y_k$} if they agree to index $k$ but not to index $k+1$.
Considering the poset in Figure~\ref{Lambda11333:fig}, we see that
\beq
\label{C:exa}
C:333\covers 133\covers 131\covers 111\covers 11
\eeq
and
\beq
\label{C':exa}
C':333\covers 133\covers 33\covers 31\covers 11 
\eeq
agree to indices 0 and 1 but not to index 2, so they diverge from index 1.

An ordering of the maximal chains of $[x,z]$
\beq
\label{<}
C_1<C_2<C_3<\cdots
\eeq
is called a {\it poset lexicographic order (PLO)\/} if it satisfies the following property.  Suppose that $C,D$ diverge from index $k$ and that  $C,C'$ agree to index $k+1$ as  do $D,D'$.  In this situation we insist that $C<D$ if and only if $C'<D'$.   This notion of chain ordering includes many of the standard ones such as EL-ordering and CL-ordering.  From now on, we assume that the maximal chains have been given a PLO.

The other ingredient in Theorem~\ref{dmt:thm} is the notion of a critical chain.   To define this, we need to look at certain intervals in a maximal chain.  Let $C$ have the form~\ree{C}.  A {\it closed interval\/} in $C$ is a subchain of the form
$$
C[y_i,y_j]: y_i\covers y_{i+1}\covers\cdots \covers y_j.
$$
Open intervals $C(y_i,y_j)$ are defined similarly.  Note that in this notation $C[y_i,y_j]$ implies $y_i\ge y_j$, while for an interval in a poset $[y_i,y_j]$ we have $y_i\le y_j$.  A {\it skipped interval (SI)\/} in $C$  is an interval $I\sbe C(z,x)$ such that $C-I\sbe B$ for some $B<C$ in the order~\ree{<}.  A {\it minimal skipped interval (MSI)\/} is an SI which is minimal with respect to containment.  For example, consider the poset in Figure~\ref{Lambda11333:fig}  again and the chain
\beq
\label{D:exa}
D:333\covers 332\covers 33\covers 31\covers 11.
\eeq
Using the PLO defined below, it will turn out that the chains $C$ in~\ree{C:exa} and $C'$ in~\ree{C':exa} satisfy $C,C'<D$.  It follows that $D[332,31]=\{332,33,31\}$ and $D[332,332]=\{332\}$ are both SIs of $D$.  The former is not an MSI since it strictly contains the latter.  The latter must be an MSI since it contains only one element.

Let $\cI(C)$ be the set of MSIs of a maximal chain $C$ in an interval $[x,z]$.
In order to define the critical chains, we need to turn $\cI(C)$ into a set of disjoint intervals.  To this end, order the intervals in $\cI(C)$ as $I_1, I_2, \ldots$ so that their left endpoints have increasing indices along $C$.  Note that there are never any ties  because there are no containments among MSIs.  Form a new set of intervals $\jint(C)$, as follows.  Let $J_1=I_1$.  Now consider $I_2'=I_2-J_1$, $I_3'=I_3-J_1$, and so forth.  Throw out any of these intervals which are not containment minimal, and let $J_2$ be the one of smallest index which is left.  Continue in this manner until there are no more intervals to consider.  Call $C$ {\it critical\/} if it is covered by $\jint(C)$ in the sense that
$$
C(z,x)=\biguplus_i J_i.
$$
In this case, the {\it critical dimension\/} of $C$ is
$$
d(C)=|\jint(C)|-1.
$$
Continuing with our example, it turns out that the chain $D$ in~\ree{D:exa} is critical since every 1-element interval in $D(333,11)$ is an MSI.  So $\cI(D)=\jint(D)$ and $d(D)=3-1=2$.
We now have everything in place to compute $\mu$.
\bth[\cite{bh:dmf}]
\label{dmt:thm}
Let $P$ be a poset and $x,y\in P$ such that $[x,y]$ is finite.  For any PLO on the maximal chains of $[x,y]$,
$$
\mu(x,y)=\sum_{C}(-1)^{d(C)},
$$
where the sum is over all critical chains $C$.
\hqed
\eth 

\subsection{Discrete Morse theory for generalized subword order}

We will now develop the ideas needed to apply Theorem~\ref{dmt:thm} to generalized subword order.
We will assume henceforth that $P_0$ is locally finite so that  the M\"obius function can be  computed for any interval of $P_0$.   In order to distinguish concepts in $P_0$ from the same concept in $P^*$, we will adjoin a subscript zero to the former.  So, for example, $\mu_0$ is the M\"obius function for $P_0$ while $\mu$ is the M\"obius function for $P^*$.

In order to use Theorem~\ref{dmt:thm} to prove Theorem~\ref{main:thm}, we henceforth assume that $P$ comes equipped with a {\it natural labeling\/}, that  is, an injection $\ell:P\ra \bbP$ such that if $a< b$ in $P$ then $\ell(a)<\ell(b)$.  We will also let $\ell(0)=0$.  (If $P$ is not countable, then we are free to use another totally ordered set in place of $\bbP$.)  When writing out examples, we will often use $x$ and $\ell(x)$ interchangeably.  But in definitions, results, and proofs we will be careful to distinguish $a\leq_0 b$ which refers to the partial order in $P_0$, and $\ell(a)\leq \ell(b)$ which refers to the total order on $\bbP$.

We need to define PLOs for both $\pzero$ and $P^*$.  We give $\pzero$ a PLO as follows: letting $\coversinp$ denote the covering relation in $\pzero$, label the covering relation $y \coversinp x$ by $\ell(x)$, and label a maximal chain $C_0$ in $\pzero$ by the sequence $L(C_0)$ of edge labels  from top to bottom.  Clearly, ordering the maximal chains of $\pzero$ by lexicographic order of these label sequences gives a PLO. 

We now describe how to extend this idea to give a PLO for $P^*$.  Among all embeddings of $u$ in $w$ there is always a {\it rightmost\/} choice $\rho$ which has the property that for any embedding $\emb$  of $u$ in $w$ and any element $a$ of $u$,  if $\emb(i)$ and $\rho(j)$ correspond to $a$ then $i\le j$.  (Here we are considering different copies of the same element of $P$ in $u$ as distinguishable.)  
For example, using the poset in Figure~\ref{Lambda:fig} as $P$, the rightmost  embedding if $1132$ in $2132333$ would be $\rho=0010132$; note that $0001132$ is not an embedding in $w$ since $1 \not\leq_0 2$ in $P_0$.
When $\emb$ is an embedding of $u$ in $w$ and $\emb(j)=0$, we say that the $j$th position 
of $w$ has been \emph{zeroed out}; in our example, the first, second and fourth positions of $w$ get zeroed out to give $\rho$.

We wish to associate a unique embedding with each cover $w\covers u$.  Note that an embedding of such a $u$ in $w$ can be obtained by replacing some element $b$ of $w$ by an element $a$ covered by $b$ in $P_0$.  If $a\neq 0$, then $\len{u}=\len{w}$ and this embedding is unique.  If $a=0$, then $\len{u}=\len{w}-1$ and there may be several embeddings.  Among these, we will always choose the rightmost.
Now given an interval $[u,w]$ and a maximal chain
\beq
\label{wordC}
C:w=v_0\covers v_1\covers \cdots \covers v_n=u,
\eeq
each cover $v_{i-1}\covers v_i$ defines an embedding of $v_i$ in $v_{i-1}$ and thus, inductively, an embedding $\emb_i$ of $v_i$ in $w$.  We label the cover with the label $\ell_i=\spn{j_i,x_i}$ where $j_i$ is the index where $\emb_{i-1}$ and $\emb_i$ differ and $x_i=\emb_i(j_i)$.  The {\it label sequence of the chain\/} is
$$
L(C)=(\ell_1,\ell_2,\ldots,\ell_n)
$$
and we also write
$$
C:\emb_0\coversl{\ell_1} \emb_1\coversl{\ell_2} \cdots \coversl{\ell_n} \emb_n.
$$
To illustrate, here are the label sequences for the chains in~\ree{C:exa}, \ree{C':exa}, and~\ree{D:exa}:
$$
\barr{l}
C:333\coversl{\spn{1,1}} 133\coversl{\spn{3,1}} 131\coversl{\spn{2,1}} 111\coversl{\spn{1,0}} 011,\\
C':333\coversl{\spn{1,1}} 133\coversl{\spn{1,0}} 033\coversl{\spn{3,1}} 031\coversl{\spn{2,1}} 011,\\ 
D:333\coversl{\spn{3,2}} 332\coversl{\spn{3,0}} 330\coversl{\spn{2,1}} 310\coversl{\spn{1,1}} 110.
\earr
$$

We can now define the desired order.  First of all, lexicographically order the labels by letting $\spn{j,x}\lexleq\spn{k,y}$ if either $j<k$, or $j=k$ and $\ell(x)<\ell(y)$  where $\ell$ is our natural labeling of $P_0$.  This induces a lexicographic ordering on label sequences.  So we let $C<D$ if and only if $L(C)\lexleq L(D)$.  It should now be clear in our running example that $C,C'<D$ as claimed above, and also that $C' < C$.  The following proposition is easily deduced from the definitions, and its proof is left the the reader.
\bpr
Let $P$ be a naturally labeled poset such that $P_0$ is locally finite.  Then the order defined on chains of $P^*$ is a PLO.\hqed
\epr

\section{Proof of the main result}\label{pmr}

Our goal for this section is to prove Theorem~\ref{main:thm}.  Although we will use classical M\"obius function techniques towards the end of the proof, the majority of the proof will involve discrete Morse theory.  Therefore, we wish to work towards classifying the MSIs and the critical chains.  
We will first prove that certain $1$-element sets are MSIs and eliminate most of the chains of $[u,w]$ as candidates for critical chains.  This will allow us to restrict our attention to maximal chains 
$C[w,u]$ where the only change from $w$ to $u$ is that a single position of $w$ is being reduced.
For a restricted class of chains, we will need to determine the M\"obius function using classical techniques, before putting everything together at the end of the section to prove Theorem~\ref{main:thm}.  There are some parallels between the start of this section and \cite[\S 5]{sv:mfc}, particularly in our Lemmas~\ref{chainspec:lem}, \ref{descent:lem} and \ref{ascent:lem}, although overall our case is significantly more involved.  

\subsection{Critical chains must have lexicographically decreasing labels}

Our goal for this subsection is to prove Corollary~\ref{lexdec:cor}, which gives the result stated in the subsection title.  Along the way, we will prove some results that will be useful both here and later.

It will be useful to permute the labels of a chain's label sequence in order to produce earlier chains and resulting MSIs.  Because of our convention of always using the rightmost embedding, such a permutation may not result in a label sequence for another chain.  But something can still be said in this situation.  First of all, we need to look more carefully at the notion of a rightmost embedding.  If $w=w(1)w(2)\wdots w(n)\in P^*$ and $a\in P$ then a {\it run of $a$'s in $w$\/} is a maximal interval of indices $[i,k]=\{i,i+1,\ldots,k\}$ such that
$$
w(i)=w(i+1)=\cdots=w(k)=a.
$$ 
For example, if $w=aaabaaccc$ then $[1,3]$ is a run of $a$'s, $[4,4]$ is a run of $b$'s, $[5,6]$ is another run of $a$'s, and $[7,9]$ is a run of $c$'s.  Now suppose $w\covers u$.  If $\len{u}=\len{w}$ then, as mentioned previously, there is a unique embedding of $u$ in $w$.  But if $\len{u}=\len{w}-1$ then $u$ must be obtained from $w$ by removing a minimal element, say $a$, of $P$.  In order to obtain the rightmost embedding, one must choose $a$ to be the element of smallest index in its run.  Taking $u=aaabaccc$ in our example, we see that this cover would be labeled as
$$
aaabaaccc\coversl{\spn{5,0}} aaab0accc.
$$

Now let $C$ be a maximal chain in $[u,w]$ with label sequence $L=L(C)$.  A permutation $L'$ of the labels of $L$ will be called {\it consistent\/} if, for every index $i$, the labels of the form $\spn{i,x}$ occur in the same order in $L'$ that they do in $L$.  To illustrate, suppose we consider the following chain in the poset of Figure~\ref{Lambda11333:fig}:
$$
C: 333 \coversl{\spn{2,1}} 313 \coversl{\spn{2,0}} 303 \coversl{\spn{1,1}} 103 \coversl{\spn{3,1}} 101.
$$
In this case, 
$$
L'=(\spn{1,1},\spn{2,1},\spn{3,1},\spn{2,0})
$$ 
would be consistent with $L$, but 
$$
L''=(\spn{1,1},\spn{2,0},\spn{3,1},\spn{2,1})
$$ 
would not.  If $L'$ is consistent with $L$, then $L'$ defines a sequence of embeddings in $w$ ending at the same embedding of $u$ as with $L$, although some of the embeddings in the sequence may not be rightmost.  These embeddings define a maximal chain $C'$ (by ignoring the zeros) which is called the {\it chain specified by $L'$\/}.  Continuing our example, we can use $L'$ to generate a sequence
$$
333 \coversl{\spn{1,1}} 133 \coversl{\spn{2,1}} 113 \coversl{\spn{3,1}} 111 \coversl{\spn{2,0}} 101
$$
and so the chain specified by $L'$ is
$$
C':333 \covers 133 \covers 113 \covers 111 \covers 11
$$
with 
$$
L(C')=(\spn{1,1},\spn{2,1},\spn{3,1},\spn{1,0}).
$$  
Although we may not have $L(C')=L'$, at the first place where the label sequences differ it must be because a position is zeroed out, and that position is further left for $L(C')$ than for $L'$.
So we have proved the following result.
\ble[Chain Specification Lemma]\label{chainspec:lem}
Let $C$ be a maximal chain in $[u,w]$ and let $L'$ be a consistent permutation of $L(C)$.  If $C'$ is the chain specified by $L'$, then $\chainlabel{C'} \lexleq L'$, and $C'$ ends at $u$. \hqed
\ele

To describe the MSIs in generalized subword order, we will need some definitions.  Let
$$
C: w=v_0 \covers v_1 \covers \cdots \covers v_n=u,
$$
or, in our alternative notation,
$$
C: \emb_0 \coversl{\spn{i_1,x_1}} \emb_1 \coversl{\spn{i_2,x_2}} \cdots \coversl{\spn{i_n,x_n}} \emb_n
$$
 be a maximal chain in $[u,w]$.  
For $1 \leq j < n$, we see that $v_j$ comes between the labels $\spn{i_j,x_j}$ and $\spn{i_{j+1},x_{j+1}}$ along $C$, and we will say that \emph{$v_j$ has the labels $\spn{i_j,x_j}$ and $\spn{i_{j+1},x_{j+1}}$}.  
With this in mind, we will say that $v_j$ is a \emph{1-descent} if $i_j > i_{j+1}$. 
This terminology is meant to remind the reader that it is the first element of a label pair which is being reduced.
 If $i_j<i_{j+1}$ we will say $v_j$ is an \emph{ascent}.  A chain $C$  is said to be {\it weakly 1-increasing\/} if $L(C)=(\spn{i_1,x_1}, \spn{i_2,x_2}, \ldots, \spn{i_n,x_n})$ satisfies $i_1 \leq i_2 \leq \cdots \leq i_n$.

Using Lemma~\ref{chainspec:lem}, we will prove the following simple but useful statement.
\ble[Descent Lemma]\label{descent:lem} A 1-descent is an MSI of one element.
\ele 

\bex
In the chain
$$
322 \coversl{\spn{2,0}} 302 \coversl{\spn{1,2}} 202,
$$
we have that $\{32\}$ is a 1-descent.  It is also the case that $\{32\}$ is an MSI because of the lexicographically earlier chain
$$
322 \coversl{\spn{1,2}} 222 \coversl{\spn{1,0}} 022.
$$
\eex

\bprf[Proof of Lemma~\ref{descent:lem}]
Suppose $C$ is a maximal chain in $[u,w]$ which
includes $v^+ \covers v \covers v^-$, and that $v$ is a 1-descent with labels $\spn{i,x}$ and $\spn{j,y}$, so $i > j$.  Let $L'$ be the label sequence $(\spn{j,y}, \spn{i,x})$ which is consistent with $L(C[v^+,v^-])$.  Form a new chain $C'$ from $C$ by replacing the interval $C[v^+,v^-]$ by the chain specified by $L'$.  By definition of $C'$  we have $C-\{v\} \subseteq C'$.  Also,  by Lemma~\ref{chainspec:lem}, 
\[
\chainlabel{C'[v^+,v^-]} \lexleq L' \lexless \chainlabel{C[v^+,v^-]},
\]
and so $C'$ lexicographically precedes $C$. Thus $\{v\}$ is an SI, which must be an MSI since it has just one element.
\eprf

We can now eliminate a large class of chains from consideration as critical chains.
\ble[Ascent Lemma]\label{ascent:lem}
If an interval $I$ of a maximal 
chain $C$ contains an ascent, then $I$ is not an MSI.
\ele

\bex
In the proof that follows, it may help the reader to consider the example with $P=\Lambda$, $w=233$, $u=2$ and 
\[
C:233 \coversl{\spn{2,1}} 213 \coversl{\spn{2,0}} 203 \coversl{\spn{3,1}} 201 \coversl{\spn{3,0}} 200.
\]
The element 23 is an ascent and $I$ can be considered to be $C(233,2)$.  Assuming $I$ is an MSI, our proof will eventually show that $C(23,2)$ is an MSI, yielding a contradiction.  Relevant to the first part of the proof is a weakly 1-increasing chain that ends at the rightmost embedding, which in this case is
\[
233 \coversl{\spn{1,0}} 033 \coversl{\spn{2,1}} 013 \coversl{\spn{2,0}} 003 \coversl{\spn{3,2}} 002.  
\] 
\eex

\bprf[Proof of Lemma~\ref{ascent:lem}]
Suppose $C$ is a maximal chain from $w$ to $u$ and an interval $I$ of $C$ contains an ascent but is an MSI.  Our eventual goal is to obtain a contradiction by showing that $I$ strictly contains an SI.  Without loss of generality, we can assume that $I = C(w,u)$.  By Lemma~\ref{descent:lem}, $I$ does not contain a 1-descent, and so $C$ is weakly 1-increasing.  Let $\emb$ be the embedding of $u$ in $w$ determined by $C$, and let $\rho$ be the rightmost embedding.  
Since we wish to first show that $\emb \neq \rho$, suppose, towards a contradiction, that $\emb = \rho$ and so $C$ ends at $\rho$.  Since $C$ has an ascent, $C$ decreases letters in at least 2 positions. Combined with the fact that $C$ is weakly 1-increasing, we see that $C$ decreases each such position using the lexicographically first chain in $\pzero$
between the starting and ending element of that position, since otherwise $I$ would not be a minimal SI.  However, the weakly 1-increasing chain to $\rho$ which uses the lexicographically first chain of $\pzero$ in each position is the lexicographically first chain in $[u,w]$, contradicting the fact that $C(w,u)$ is an MSI.  We conclude that $\emb$ is not rightmost.

For an embedding $\zeta$ of $u$ in $w$ and all indices $j$ with $1 \leq j \leq \len{w}$, define
$\zeros_\zeta(j)$ to be the number of indices $i$ with $i \leq j$ such that $\zeta(i)=0$.  Because $\rho$ is rightmost, we have $\zeros_\rho(j) \geq \zeros_\emb(j)$ for all $j$ and we have equality when $j = \len{w}$.  Since $\rho \neq \emb$, there must be some first index $a$ where $\zeros_\rho(a) > \zeros_\emb(a)$.  Since we have equality when $j = \len{w}$, there must be a first index $c > a$ such that $\zeros_\rho(c) = \zeros_\emb(c)$.  Note that by the choice of $a$ and $c$, we must have $\rho(a)=0$, $\emb(a) \neq 0$, $\rho(c)\neq 0$ and $\emb(c)=0$.  Letting $\rho(c)=x$ and $\emb(a) =y$, we have the situation in the following diagram, with $b$ to be explained next:
\[
\begin{array}{lcccccccc}
w:  & \quad & \wdots & w(a) & \wdots & w(b) & \wdots & w(c) & \wdots \\
\rho: & & \wdots & 0 & \wdots &  & \wdots & x & \wdots \\
\emb: & & \wdots & y & \wdots & x & 0 \wdots0 & 0 & \wdots \\
\end{array}.
\]
Since $\rho$ and $\emb$ are embeddings of the same word $u$ and by the choice of $a$ and $c$, there must be an index $b$ with $a \leq b < c$ such that $\emb(b) = \rho(c) = x$ and $\emb(j)=0$ for $b < j \leq c$.  In other words, $x$ represents a given letter of $u$ in two different positions in $\emb$ and $\rho$.  Note that $w(b), w(c) \geq_0 x$ 
since $\rho$ and $\emb$ are embeddings.

We next show that in going from $w$ to $u$ along $C$, the conditions on $C$ imply that only those letters of $w$ in positions $i$ with $b < i \leq c$ are decreased.  To do this, define $p, q \in C$ to be the elements such that $C[q,p]$ contains all the labels of the form $\spn{i,*}$ for $b < i \leq c$.  Note that since $C$ is weakly 1-increasing, we must have
 $\zeta(b)=x$ and $\zeta(c)=w(c) \geq_0 x$, 
where $\zeta$ is the embedding corresponding to $q$ on $C$.
Define a sequence of labels $L'$ by
first zeroing out position $b$ in $\zeta$, then zeroing out positions $i$ for $b < i< c$ from left to right but otherwise in any way, and finally decreasing $\zeta(c)$ to $x$.  
Form a new chain $C'$ by replacing the interval $C[q,p]$ of $C$ by the chain specified by $L'$.  By Lemma~\ref{chainspec:lem} and the construction of $L'$, we have $\chainlabel{C'[q,p]} \lexleq L' \lexless \chainlabel{C[q,p]}$ since $L'$ zeroes out position $b$ but $L(C[q,p])$ does not change position $b$.  Thus $C'$ lexicographically precedes $C$.  Hence $C(q,p)$ is an SI of $C$, and since $C(w,u)$ is an MSI, this forces $q=w$ and $p=u$.  Also, $w(b) = q(b) = x$.  

To obtain a contradiction to our original assumption that $I=C(w,u)$ is an MSI, we finish by showing that $C$ strictly contains an SI.  Since $C(w,u) = C(q,p)$ we have reduced to the case where $C$ reduces neither positions $i\leq b$ nor positions $i>c$.  Since $C$ has an ascent, there must be a position $i$ with $b < i < c$ which $C$ reduces.  By definition of $b$, $C$ must zero out all such positions as well as position $c$.  Since $C$ is weakly 1-increasing, it must pass through the word 
\[
w' = (w(1), \ldots, w(b-1), x, w(c), w(c+1), \ldots, w(\len{w})).
\]
Because there are positions $i$ for $b<i<c$ to be zeroed out, $w' \neq w$.  Note that $u$ is obtained from $w'$ on $C$ by decreasing $w(c)$ to 0.  Define a label sequence $L'$ by first zeroing out position $b$ (containing the letter $x$) in $w'$ and then decreasing $w(c)$ to $x$.  Form a new chain $C'$ by replacing the interval $C[w',u]$ by the chain specified by $L'$.  By Lemma~\ref{chainspec:lem} and the construction of $L'$, $\chainlabel{C'[w',u]} \lexleq L' \lexless \chainlabel{C[w',u]}$.  Thus $C'$ lexicographically precedes $C$, and $C(w',u)$ is an SI in $C$. This yields the desired contradiction since $w' < w$ and $C(w,u)$ is an MSI.  
\eprf

\bco\label{lexdec:cor}
The labels on any critical chain must be lexicographically decreasing.
\eco

\bprf
A critical chain must be covered by MSIs and so, by Lemma~\ref{ascent:lem}, cannot contain an ascent.  Therefore, if $v^+ \covers v \covers v^-$ appears on a critical chain $C$ and $v$ has labels $\spn{i,x}$ and $\spn{j,y}$, we know $i \geq j$.  If $i >j$, we have $\spn{i,x} \lexgreater \spn{j,y}$.  If $i=j$, since $v \covers v^-$, we must have $x \coversinp y$ and hence $\ell(x) > \ell(y)$, where $\ell$ is our natural labeling of $P_0$.  Again, we get $\spn{i,x} \lexgreater \spn{j,y}$, as required.
\eprf

\subsection{Classification of the most important MSIs}

As a consequence of Lemma~\ref{descent:lem} and Corollary~\ref{lexdec:cor}, we can restrict our attention to MSIs of a very special form.  Consider a potentially critical chain $C$ from $w$ to an embedding $\emb$, and let $k$ denote the largest index such that $w(k) \neq \emb(k)$.  Since the labels must be lexicographically decreasing, $C$ must start by reducing $w(k)$ to $\emb(k)$.  Then $C$ must reduce $w(j)$ to $\emb(j)$ where $j < k$ is as large as possible with $w(j)\neq \emb(j)$, and so on from right to left.  A key observation is that, by Lemma~\ref{descent:lem}, in moving from decreasing letters in position $k$ to decreasing letters in position $j$ to the left of position $k$, we will create a 1-descent in the label sequence and hence a single-element MSI.  Thus the only MSIs left to determine on potentially critical chains are those of the form $C(w,u)$, where the only change from $w$ to $u$ is that a single position has been reduced, including the mathematically trickier possibility of it being zeroed out.  The proposition below classifies these MSIs.  

To state the proposition, we will need two natural ideas.  For the first, 
suppose $C$ is a maximal chain in $[u,w]$ and that the only difference between $w$ and an embedding $\emb$ of $u$ is that the $j$th position has been decreased.  Then there is an obvious bijection from the maximal chains of the interval $[u,w]$ in $\pstar$ that end at $\emb$ to the maximal chains of the interval $[\emb(j), w(j)]$ in $\pzero$, namely the isomorphism that sends each embedding $\zeta$ to $\zeta(j)$.  Throughout, let $\czero$ denote the image of $C$ under this bijection.
  
The second idea needed to state the proposition is given by the following lemma where, similar to our earlier convention, the condition $w(j-1) \leq_0 w(j)$ is considered false when $j=1$ 
since $w(0)$ does not exist.

\ble\label{rightmost:lem}
Suppose $\emb$ is an embedding of $u$ in $w$ and that $w$ and $\emb$ differ only in their $j$th position.  Then $\emb$ is not the rightmost embedding of $u$ in $w$ if and only if $\emb(j)=0$ and $w(j-1) \leq_0 w(j)$.  
\ele

\bprf
Suppose $\emb$ is not the rightmost embedding $\rho$ of $u$.  Then $\emb(j)=0$ since otherwise there is only one embedding of $u$.  Since $\emb$ has exactly one 0, the fact that $\rho$ is further right than $\emb$ implies that $\rho(j)\neq 0$ and $\rho(i)=0$ for some $i<j$. It follows that $w(j) \geq_0 \rho(j) =\emb(j-1)=w(j-1)$.

For the converse, suppose $\emb(j)=0$ and $w(j-1) \leq_0 w(j)$.  To show that $\emb$ is not rightmost, define $\emb'$ as the result of swapping $\emb(j-1)$ and $\emb(j)$, i.e., $\emb'(j-1)=\emb(j)=0$ and $\emb'(j)=\emb(j-1)$ and $\emb'(i)=\emb(i)$ for $i \neq j, j-1$.  To show that $\emb'$ is an embedding of $u$ in $w$, we need only show that $\emb'(j) \leq_0 w(j)$, which follows from $\emb'(j) = \emb(j-1) \leq_0 w(j-1) \leq_0 w(j)$.  So $\emb'$ is an embedding further to the right than $\emb$.
\eprf

For a maximal chain $C$ in an interval $[x,y]$ of a poset, define an SI to be \emph{proper} if it is strictly contained in $C(y,x)$.  Also, if $\emb$ is the embedding of $u$ in $w$ determined by a chain $C$ we will sometimes write $C(w,\emb)$ for $C(w,u)$.

\bpr\label{msi:pr}
Suppose an embedding $\emb$ differs from $w$ only in its $j$th position.  
If $C$ is a maximal
chain from $w$ to $\emb$, then $C(w,\emb)$ is an MSI if and only if either
\begin{enumerate}
\item[(a)] $\emb$ is rightmost and $\czero(w(j),\emb(j))$ is an MSI of $\czero$ in $\pzero$, or
\item[(b)] $\emb$ is not rightmost and $\czero(w(j),\emb(j))$ contains neither an element $x \geq_0 w(j-1)$ nor a proper SI of $\czero$.
\end{enumerate}
\epr

This proposition tells us that if $\emb$ is rightmost, for which we have an easy test by Lemma~\ref{rightmost:lem}, then MSIs of the form $C(w,\emb)$ are in perfect correspondence with those of the form $\czero(w(j),\emb(j))$ in $\pzero$.  If $\emb$ is not rightmost, then things are more complicated.  It  may be easier to digest Condition (b) if we can divide it into two cases: $\czero(w(j),\emb(j))$ does not contain a proper SI if and only if either it is an MSI itself, or it is not an MSI itself and it contains no MSIs.  In the latter case, $\czero$ must  be the lexicographically first chain in $\pzero$ from $w(j)$ to $\emb(j)$.  So we can restate (b) as
\textit{
\begin{enumerate}
\item[(b$'$)] $\emb$ is not rightmost and $\czero(w(j),\emb(j))$ does not contain an element $x \geq_0 w(j-1)$ and is either an MSI or is lexicographically first.
\end{enumerate}
}

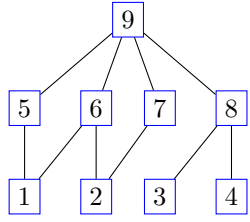
\begin{figure}
\begin{center}
\begin{tikzpicture}
\tikzstyle{elt}=[rectangle, draw=boxes]
\matrix{
&[15pt]&\node(9)[elt]{$9$};&&[15pt]\\[20pt]
\node(5)[elt]{$5$};&\node(6)[elt]{$6$};&&\node(7)[elt]{$7$};&\node(8)[elt]{$8$};\\[20pt]
\node(1)[elt]{$1$};&\node(2)[elt]{$2$};&&\node(3)[elt]{$3$};&\node(4)[elt]{$4$};\\
};
\draw (9)--(5);
\draw (9)--(6);
\draw (9)--(7);
\draw (9)--(8);
\draw (5)--(1);
\draw (6)--(1);
\draw (6)--(2);
\draw (7)--(2);
\draw (8)--(3);
\draw (8)--(4);
\end{tikzpicture}
\caption{The poset $P$ of Example~\ref{MSI:ex}}
\label{MSI_example:fig}
\end{center}
\end{figure}

\begin{figure}
\begin{center}
\begin{tikzpicture}
\tikzstyle{elt}=[rectangle, draw=boxes]
\matrix{
&[20pt]&[20pt]&[20pt]\node(29)[elt]{$29$};&[20pt]&[20pt]\\[30pt]
\node(9)[elt]{$9$};&&\node(25)[elt]{$25$};&\node(26)[elt]{$26$};&\node(27)[elt]{$27$};&\node(28)[elt]{$28$};\\[30pt]
\node(6)[elt]{$6$};&\node(7)[elt]{$7$};
      &\node(21)[elt]{$21$};&\node(22)[elt]{$22$};&\node(23)[elt]{$23$};&\node(24)[elt]{$24$};\\[30pt]
&&&\node(2)[elt]{$2$};&&\\
};
\draw (29)--(9);
\draw (29)--(25);
\draw (29)--(26);
\draw (29)--(27);
\draw (29)--(28);
   \draw (9)--(6);
   \draw (9)--(7);
   \draw (25)--(21);
   \draw (26)--(6);
   \draw (26)--(21);
   \draw (26)--(22);
   \draw (27)--(7);
   \draw (27)--(22);
   \draw (28)--(23);
   \draw (28)--(24);
\draw (2)--(6);
\draw (2)--(7);
\draw (2)--(21);
\draw (2)--(22);
\draw (2)--(23);
\draw (2)--(24);
\end{tikzpicture}
\caption{The interval $[2,29]$ with $P$ from Figure~\ref{MSI_example:fig}}
\label{special_case_example:fig}
\end{center}
\end{figure}
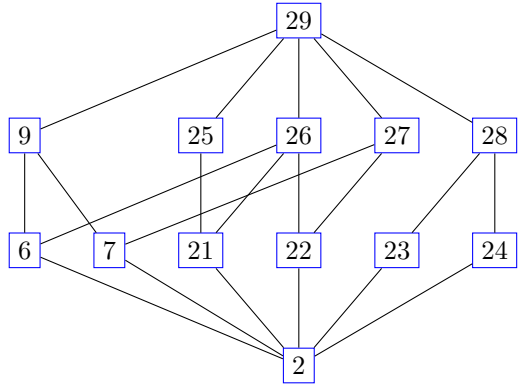

\bex\label{MSI:ex}
Consider the poset $P$ of Figure~\ref{MSI_example:fig}.  Of course, $\pzero$ can be obtained from $P$ by adding $0$ as a bottom element. To gain some intuition for Proposition~\ref{msi:pr}, look at  the interval in $P^*$ displayed in Figure~\ref{special_case_example:fig}.   Let us first consider the rightmost (and only) embedding $\emb=21$ of $u=21$ in $w=29$.  The single-element MSI $\{26\}$ of $29 \covers 26 \covers 21$ in $\pstar$ corresponds exactly to the single-element MSI $\{6\}$ of $9 \coversinp 6 \coversinp 1$ in $\pzero$, consistent with Condition~(a).  As a more complicated example, let us use the proposition to determine those maximal chains $C$ for which the entire open interval $C(29,20)$ is a single MSI.  This embedding $20$ is not rightmost, as can be seen directly or checked by Lemma~\ref{rightmost:lem}.
By Condition~(b), any maximal chains $C$ from $29$ to $20$ through $26$ or $27$ will not contain an MSI of cardinality 2 since $\czero(9,0)$ will contain an element $x=6$ or $7$, and hence an element $x \geq_0 2$.  When $C=29 \covers 28 \covers 24 \covers 20$, we see that $\czero(9,0)$ contains a proper SI, namely \{4\}, so $C(29,20)$ again violates Condition~(b).  
The remaining two maximal chains of $[2,29]$ which end at 20 satisfy Condition~(b); they are
$29\covers 25\covers 21\covers 20$ and $29\covers 28\covers 23\covers 20$, and yield the two-element MSIs $\{25,21\}$ and $\{28,23\}$ respectively.  In the context of Condition~(b$'$), we note that $9 \coversinp 5 \coversinp 1\coversinp 0$ is the lexicographically first maximal chain in $\pzero$ and corresponds to the first of our 2-element MSIs.  Regarding the second two-element MSI, $9 \coversinp 8 \coversinp 3\coversinp 0$ is a chain in $\pzero$ with no proper MSI and so corresponds to the other case in that clause of (b$'$).
It is worth noting that there are 2 other critical chains in $[2,29]$ but they each contain two single-element MSIs.  These 2 additional critical chains are
\[
\begin{array}{c}
 29 \covers 26 \covers 21 \covers 20, \\
  29 \covers 27 \covers 22 \covers 02.
 \end{array}
 \]
\eex

\begin{proof}[Proof of Proposition~\ref{msi:pr}]
Let $u \in \poset^*$ be the word corresponding to the embedding $\emb$.  Note that $C(w,\emb)$ is an MSI if and only if $C$ is not the lexicographically first chain in $[u,w]$ and $C$ has no proper skipped intervals; the analogous statement holds for $\czero(w(j), \emb(j))$.   We will first prove both directions of the proposition for the case when $\emb$ is rightmost, and then prove both directions for the case when $\emb$ is not rightmost.

Assume first that $\emb$ is rightmost because $\emb(j) \neq 0$.  
Then $\emb$ is the only embedding of $u$ in $w$, and so there is an isomorphism from $[u,w]$ in $\pstar$ to $[\emb(j),w(j)]$ in $\pzero$ which sends $v \in [u,w]$ to $\zeta(j)$, where $\zeta$ is the embedding corresponding to $v$.  Since this isomorphism preserves the lexicographic ordering of maximal chains, the result follows in this subcase.

Now suppose $\emb$ is rightmost and $\emb(j)=0$.  Since we are still in the case that $\emb$ is rightmost, we must have $w(j-1) \not\leq_0 w(j)$ by Lemma~\ref{rightmost:lem}.  By definition, $C$ leaves the elements $w(1), \ldots, w(j-1)$ of $w$ fixed.  We claim that any chain $B \lexless C$ from $w$ to $u$ must also leave these elements fixed.  Indeed, suppose $B$ decreases $w$ to the embedding $\emb'$ of $u$ in $w$ and, along the way, $B$ decreases $w(k)$ with $k<j$.  Since $B$ and $C$ both end at $u$, the first $j-1$ nonzero elements of $\emb'$ must be the letters $w(1), \ldots, w(j-1)$ in that order. So if $k$ is the leftmost position modified by $B$, we see that $B$ must decrease $w(k)$ to 0.  Moreover, since $u$ has just one less letter than $w$, $B$ decreases exactly one letter to 0.  As a result, we deduce that the $(j-1)$st nonzero elements of $\emb$ and $\emb'$ are $w(j-1)$ and $\emb'(j)$ respectively.  Thus $w(j-1) = \emb'(j) \leq_0 w(j)$, which is a contradiction.  

We have determined that $B$ fixes $w(1), \ldots, w(j-1)$.  Since $B \lexless C$, just after $B$ and $C$ diverge, the letters in the $j$th place, $b$ in $B$ and $c$ in $C$, must satisfy $\ell(b)<\ell(c)$.  Thus there is a chain $B_0$ lexicographically earlier than $\czero$ in $\pzero$.  The converse is also true: the existence of such a $B_0$ implies the existence of $B$ with $B \lexless C$.  By the observation in the first paragraph, we finish this subcase by showing that $C$ has a proper skipped interval if and only if $\czero$ does.  The reverse implication is clear.  So suppose $C$ has a proper skipped interval $I$ where $I=C(w',u')$ and $\emb''$ is the corresponding embedding of $u'$ in $w$.  If $\emb''(j) \neq 0$, then $I$ corresponds to an SI  in $\pzero$ by the idea in the second paragraph.  If $\emb''(j) = 0$ then this forces $u'=u$.  Thus, since $I$ is proper, $I$ does not contain the element $v$ of $C$ covered by $w$.  Therefore, the chain $B$ of $[u,w]$ with $B \lexless C$ that causes $I$ to be an SI satisfies $(C - I) \subseteq B$ with $v \in B$.  
Now we can construct a proper SI of $\czero$ by considering the following chain of $[0,w(j)]$ in $P_0$: start at $w(j)$, follow the $j$th letter along $B$ through $v(j)$ and all the way to $b$ (where $b$ plays the same role as it did at the start of this paragraph), and then continue along any maximal chain from $b$ to $0$ in $\pzero$.  Since the resulting chain and $\czero$ both contain $v(j)$, and since $\ell(b) < \ell(c)$, we have that $\czero$ contains a proper SI.  
We conclude that when $\emb$ is rightmost, $C(w,\emb)$ is an MSI in $\pstar$ if and only if $\czero(w(j),\emb(j))$ is an MSI of $\czero$ in $\pzero$.

For the remainder of the proof, assume $\emb$ is not rightmost so, by Lemma~\ref{rightmost:lem}, $\emb(j)=0$ and $w(j-1)\leq_0 w(j)$.  To go in the forward direction, suppose $C(w,\emb)$ is an MSI.  Thus $\czero$ cannot contain a proper SI since this would create an SI in $C$ that was smaller than $C(w,\emb)$.  If an element $x$ with $w(j-1)\leq_0 x <_0 w(j)$ as in the statement of the proposition exists, then we can use it to create a smaller SI in $C$ as follows.  
Construct a chain $B$ by first following $C$ down to the embedding $\emb'$ which has $x$ as its $j$th letter.
Then, in the terminology of Lemma~\ref{chainspec:lem}, let $B$ follow the chain specified by the label sequence obtained in the following way: reduce $\emb'(j-1)$ to 0 and then reduce $\emb'(j)=x$ to $w(j-1)$. 
By construction, $B \lexless C$ and so creates an SI of $C$.  Furthermore, since $x<_0 w(j)$, this SI is proper which contradicts our assumption.  We conclude that if $C(w,\emb)$ is an MSI, then (b) is satisfied when $\emb$ is not rightmost.

For the reverse implication, suppose that (b) holds.  Since $\emb(j)=0$ and $w(j-1) \leq_0 w(j)$, we can create a chain $B \lexless C$ intersecting $C$ only at $w$ and $u$ in a manner similar to that in the previous paragraph: 
let $B$ be the chain specified by the label sequence obtained by first decreasing $w(j-1)$ to 0 and then decreasing $w(j)$ to $w(j-1)$.
Thus $C(w,\emb)$ is an SI.  To show that $C(w,\emb)$ is an MSI, assume to the contrary that $C$ contains a proper skipped interval $I$.  Let $B$ be the chain that causes $I$ to be an MSI, meaning $B \lexless C$ with $C-I \sbe B$.  If $B$ diverges from $C$ by decreasing an element in the $j$th position, then we can use $B$ to create a proper SI of $\czero$ as in the fourth paragraph of this proof.  So consider what happens if $B$ diverges from $C$ by decreasing an element in the $k$th position with $k < j$.  This is the situation of the third paragraph, where we deduced that once $B$ arrives at the embedding $\emb'$ corresponding to the word $u$, it has reduced $w(j)$ to $\emb'(j)=w(j-1)$.  But now consider the element $x$ in the $j$th position of the word $v$ from which $B$ and $C$ diverge.  We must have that $x <_0 w(j)$ since otherwise $I$ would not be proper: once $B$ and $C$ diverge, they do not meet again until they have both decreased exactly one letter to 0, which first happens when they arrive at $u$.  Thus $x \in C_0(w(j), \eta(j))$.  Since $v$ is on $B$, and $B$ reduces $w(j)$ ultimately to $w(j-1)$, we have $x = v(j) \geq_0 w(j-1)$, a final contradiction.
\end{proof}

As we discussed, the only MSIs on potentially critical chains left to determine were those of Proposition~\ref{msi:pr}.  Therefore, we might hope that we would now be ready to classify all critical chains in $[u,w]$ and hence prove Theorem~\ref{main:thm}.  However, this is not the case, which is one of the subtler issues of our proof.  As we will see, Proposition~\ref{msi:pr} will be crucial to classifying the critical chains, but we will first take care of a very special but tricky case. 

\subsection{A class of intervals requiring special treatment}

For any $u, w \in \pstar$, any critical chain from $w$ down to $u$ ends at a unique embedding $\emb$.  With Theorem~\ref{dmt:thm} in mind, let us write $\muembpstar{\emb}{w}$ for the contribution to $\mupstar(u,w)$ from the critical chains that end at $\emb$.  Theorem~\ref{dmt:thm} tells us that 
\beq\label{sumoveremb:eq}
\mupstar(u,w) = \sum_\emb  \muembpstar{\emb}{w},
\eeq
where the sum is over all embeddings $\emb$ of $u$ in $w$.  Therefore, to prove Theorem~\ref{main:thm}, we need to determine $\muembpstar{\emb}{w}$.

As we will show in the proof  of Theorem~\ref{main:thm}, the case we must consider separately involves two-letter words and is as follows: for $a, b \in \pzero$ with $a\leq_0 b$, let $w=ab$ and $\emb=a0$.  We treat this case separately because it is the one where, along with discrete Morse theory, we will use classical M\"obius function techniques.  We saw an example of this case in Example~\ref{MSI:ex} when we considered $C(29,20)$; the interval $[2,29]$ in $\pstar$ is shown in Figure~\ref{special_case_example:fig}.  Our primary goal in this subsection is to find a formula for $\muembpstar{a0}{ab}$, which we will do in Corollary~\ref{special_case:co} using the next two lemmas.

For an embedding $\emb$ of $u$ in $w$, define a subposet $[\emb, w]$ of $[u,w]$ by
\beq
\label{embinterval:eq}
\barr{ll}
[\emb, w]  =  \{v \in \pstar : &\mbox{there exists an expansion $\zeta$ of $v$ with $\len{\zeta}=\len{w}$}\\
 &\mbox{and  $\emb(j) \leq_0 \zeta(j) \leq_0 w(j)$  for all $j$}\}.
\earr
\eeq
For example, with $P$ as in Figure~\ref{MSI_example:fig}, $[20,26]$ is a poset on $\{2,21, 22, 26\}$.
A word of warning is in order here: it is generally not the case that $\muembpstar{\emb}{w}$
is the same as $\mu_{[\emb,w]}(\emb,w)$.   
This is because $[u,w]$ contains maximal chains  that $[\emb,w]$ does not, so SIs in $[u,w]$ might not be SIs in $[\emb,w]$.  Continuing with our example, $[2,26]$ includes the maximal chain $26 \covers 6 \covers 2$ but $6 \not\in [20,26]$.  As a result, $\{21\}$ is an SI of the chain $26\covers 21\covers 2$ in $[2,26]$ but not in $[20,26]$.

For $a\neq0$, critical chains contributing to $\mupstar(a,ab)$ either contribute to 
\linebreak
$\muembpstar{0a}{ab}$ or to $\muembpstar{a0}{ab}$, but not to both.  Therefore, to obtain the information we want about $\muembpstar{a0}{ab}$, we will first consider $\muembpstar{0a}{ab}$. 
This gives an example of a situation where $\muembpstar{\emb}{w}$ and $\mu_{[\emb,w]}(\emb,w)$ are in fact equal.

\ble\label{cartesian:lem}
Suppose $a \leq_0 b$ in $\pzero$.  Then 
\[
\muembpstar{0a}{ab} = \mup(0,a)\mup(a,b) = \mu_{[0a,ab]}(0a,ab).
\]
\ele

\bprf
If $a=0$ then the result is clear, so assume $a \neq 0$.  
We start with the first equality, and use discrete Morse theory.  Let $C$ be a critical chain from $ab$ to $0a$ in $\pstar$. 

 Let us first suppose that $a<_0 b$. By Corollary~\ref{lexdec:cor}, $C$ must contain $aa$.  By Lemma~\ref{descent:lem}, $\{aa\}$ is a single-element MSI.  Thus $C$ consists of any critical chain $C_1$ of $[aa,ab]$, followed by the MSI $\{aa\}$, followed by any critical chain $C_2$ of $[0a,aa]$.  The respective embeddings $aa$ in $ab$, and $0a$ in $aa$, are rightmost.  Therefore, by Proposition~\ref{msi:pr}, the MSIs, and hence the critical chains, of $[aa,ab]$ (respectively $[0a,aa]$) in $\pstar$ are in bijection with those in $[a,b]$ (resp.\ $[0,a]$) in $\pzero$, and the number of $\jint$-intervals required to cover the critical chains will be equal.  It follows that
$$
|\jint(C)|=|\jint(C_1)|+|\jint(C_2)|+1.
$$
As a result, summing over critical chains $C_1$ and $C_2$ of $[aa,ab]$ and $[0a,aa]$ respectively, Theorem~\ref{dmt:thm} gives
\bea
\muembpstar{0a}{ab} & = & \sum_{C_1, C_2} (-1)^{d(C)} \\
& = & \sum_{C_1, C_2} (-1)^{|\jint(C_1)|+|\jint(C_2)|} \\
& = & \left( \sum_{C_1} (-1)^{d(C_1)} \right) \left(\sum_{C_2} (-1)^{d(C_2)} \right)\\
& = & \mup(a,b) \mup(0,a).
\eea

If $a=b$, then $C$ will consist of just $C_2$.  By a similar argument to that in the previous paragraph, we obtain
$$
\muembpstar{0a}{ab} = \muembpstar{0a}{aa}=\mup(0,a) =
\mup(0,a)\mup(a,a)=\mup(0,a) \mup(a,b).
$$ 

While we could prove that $\muembpstar{0a}{ab} = \mu_{[0a,ab]}(0a,ab)$ using discrete Morse theory, we will instead prove the second equality using more classical techniques.  In fact, we will prove the stronger result that for $a\leq_0 b$, $[0a,ab]$ is isomorphic to the Cartesian product $[0,a] \times [a,b]$ of intervals of $\pzero$.  Let $c,d,e \in \pzero$.  Consider the map $f: [0a,ab] \to [0,a] \times [a,b]$ defined by 
\bea
f(cd) & = & (c,d) \mbox{\ \ if neither $c$ nor $d$ is 0, while} \\
f(d) & = & (0,d).  
\eea
Clearly, $f$ is injective.  Let $(c,d) \in [0,a] \times [a,b]$.  If $c,d$ are nonzero, then $cd \in [0a,ab]$.  Combined with the fact that $d \geq_0 a >_0 0$, we get that $f$ is surjective.  Observe that $f$ is order-preserving on words of the same length. 
For words of different length, suppose $e\le cd$ in $[0a,ab]$.  So $e\le_0 c$ or $e\le_0 d$.  In the former case we have, using~\ree{embinterval:eq} with the embedding $\zeta=cd$, that $e\le_0 c\le_0 a \le_0 d$.  So in either case $e\le_0 d$,
yielding $f(e) = (0,e) \leq (c,d) = f(cd)$.  Showing that $f^{-1}$ is order-preserving is similarly easy.  The result now follows by the product theorem for the M\"obius function \cite[Proposition 3.8.2]{sta:ec1}.
\eprf

To compute our quantity of interest, $\muembpstar{a0}{ab}$, we will need one more lemma which is the following general result.  Recall that an \emph{upper order ideal} or \emph{dual order ideal} or \emph{filter} of a poset $Q$ is a subposet $U$ such that if $x, y \in Q$ with $x \leq y$, then $x \in U$ implies $y \in U$.  For any poset $Q$, let $\widehat{Q}$ denote $Q$ with a bottom element $\widehat{0}$ and a top element $\widehat{1}$ adjoined.  In this case, we can abbreviate $\mu_{\widehat{Q}}(\widehat{0}, \widehat{1})$ as $\mu(\widehat{Q})$.

\ble\label{incexc:lem}
Consider a finite poset $Q$ of the form $U \cup V$, where $U$ and $V$ are upper order ideals of $Q$. 
Then 
\beq\label{incexc:eq}
\mu(\widehat{Q}) = \mu(\widehat{U}) + \mu(\widehat{V}) - \mu(\widehat{U \cap V}). 
\eeq
\ele

\bprf
We will use the expression for the M\"obius function as an alternating sum of chain counts \cite[Proposition~3.8.5]{sta:ec1}: 
\beq\label{mobius:eq}
\mu(\widehat{Q}) =  -1 + c_0 -c_1 + c_2 - c_3 + \cdots,
\eeq
where $c_i$ is the number of chains in $Q$ of length $i$ (i.e., containing $i+1$ elements).  Consider an arbitrary chain $\arbchain$ in $Q$ and let $q$ denote the smallest element of $\arbchain$.   Obviously, $q$ will be an element of exactly one of $U \setminus V$, $V \setminus U$ or $U \cap V$.  If $q \in U \setminus V$ then, because $U$ is an upper order ideal, $\arbchain$ contributes to $\mu(\widehat{U})$ with a sign determined by \eqref{mobius:eq}.  An analogous statement holds if $q\in V\setminus U$.  If $q \in U \cap V$, then $\arbchain$ contributes to both $\mu(\widehat{U})$ and $\mu(\widehat{V})$, and to $\mu(\widehat{U \cap V})$.  Now \eqref{incexc:eq} follows.
\eprf

\bco\label{special_case:co}
Suppose $a \leq_0 b$ in $\pzero$.  Then 
\[
\muembpstar{a0}{ab} = \left\{
\begin{array}{ll}
\mup(0,b)+1 & \mbox{if $a=b>_0 0$,} \\
\mup(0,b) & \mbox{if $a<_0 b$ or $a=b=0$}.
\end{array}
\right.
\]
\eco
\bprf
If $a=0$ then the result is clear, so assume $a \neq 0$.  Let us apply Lemma~\ref{incexc:lem} with $\widehat{Q}=[a,ab]$, meaning that $Q$ is the open interval $(a,ab)$.  Let $U$ be $[0a,ab]$, as defined in \eqref{embinterval:eq}, with the top and bottom elements removed.  Thus
$$
\mu(\widehat{U}) = \mu_{[0a,ab]}(0a,ab) = \muembpstar{0a}{ab}
$$ 
by Lemma~\ref{cartesian:lem}. 

To use Lemma~\ref{incexc:lem}, we must show that $U$ is an upper order ideal in $Q$.  So suppose $x\in U$ and $y\in Q$ with $y > x$.  Note that $1\le|x|\le|y|\le2$.  To show that $y\in U$, we must consider three cases depending on the cardinalities of $x$ and $y$.  We will  write out the case $|x|=1$ and $|y|=2$ as the others are similar and simpler.  So suppose $x=c$ and $y=de$ for some $c,d,e\in P$.  Then the expansion of $c$ which satisfies~\ree{embinterval:eq} is $\zeta=0c$.  It follows that $a<_0 c \le_0 b$, where the strict inequality comes from the fact that $a\not\in U$.  Since there is only one embedding of $de$ in $ab$, we must show that it satisfies~\ree{embinterval:eq} to conclude that $de\in U$.  Since $de\in Q$, we have $de \neq ab$ with $d\le_0 a$ and $e\le_0 b$, establishing half of the needed inequalities.  Furthermore, from what we have shown, $d\le_0 a <_0 c$ and so $c0$ is not an embedding in $de$.  But $c<de$ and so the only other possibility is that $\zeta=0c$ is an embedding in $de$.  Thus $a<_0 c\le_0 e$, and of course $0\le_0 d$, establishing the remaining two inequalities.
 
Let $V$ be $[a0,ab]$ with the top and bottom elements removed.  
By definition of $[a0,ab]$, any element of $V$ must take the form $ad$ where $0 <_0 d <_0 b$, and it is easy to see that $V$ is an upper order ideal of $Q$.
We see that $U \cup V=Q$, and 
$U \cap V = [aa,ab]$ with the top element removed.  If $a=b$ then $U \cap V = \emptyset$ and so $\mu(\widehat{U \cap V}) = -1$.  If $a<_0 b$, then since $U \cap V$ has a single minimal element $aa$, we get that $\mu(\widehat{U \cap V}) = 0$.  
Applying Lemma~\ref{incexc:lem} we get
\beq\label{applylemma:eq}
\mupstar(a,ab) = \left\{
\begin{array}{ll}
\muembpstar{0a}{ab} + \mu_{[a0,ab]}(a0,ab)+1 & \mbox{if $a=b$,}\\
\muembpstar{0a}{ab} + \mu_{[a0,ab]}(a0,ab) & \mbox{if $a<_0 b$}.
\end{array}
\right.
\eeq
We know that critical chains contributing to $\mupstar(a,ab)$ either contribute to 
\linebreak
$\muembpstar{0a}{ab}$ or $\muembpstar{a0}{ab}$, but not both, and so 
\beq\label{splitmu:eq}
\mupstar(a,ab) = \muembpstar{0a}{ab} + \muembpstar{a0}{ab}. 
\eeq
Combining \eqref{applylemma:eq} and \eqref{splitmu:eq} yields
\beq
\label{2mus}
\muembpstar{a0}{ab} = \left\{
\begin{array}{ll}
\mu_{[a0,ab]}(a0,ab)+1 & \mbox{if $a=b$,} \\
\mu_{[a0,ab]}(a0,ab) & \mbox{if $a<_0 b$}.
\end{array}
\right.
\eeq
The map that sends each $ad$ in $[a0,ab]$ to $d$ is an isomorphism from $[a0,ab]$ to the closed interval $[0,b]$ of $\pzero$.  Thus $\mu_{[a0,ab]}(a0,ab)=\mup(0,b)$.  The result now follows.
\eprf

\begin{exa}\label{special_case_example:ex}
Figure~\ref{special_case_example:fig} shows the interval $[2,29]$ where $P$ is as shown in Figure~\ref{MSI_example:fig}.  
We see that $U=\{6,7,9,22,26,27\}$, $V=\{21,22,23,24,25,26,27,28\}$ and so $U \cap V = \{22,26,27\}$.  

Notice that for calculating $\muembpstar{20}{26}$,  the only critical chain  is $C: 26 \lra 21 \lra 2$ whose MSI is created by the lexicographically earlier chain $B: 26\lra 6\lra 2$.  When we consider $[20,26]$, the element $6$ is no longer present, 
but $22$ is now included. So $C$ becomes the lexicographically first chain in the interval, but the later chain, $D:26\lra 22\lra 2$ is now under consideration, when it was not before, and is critical.
So even though the critical chain is different in the two cases, there number of critical chains is the same, consistent with equation~\ree{2mus}.  We should note that more complicated examples can be constructed where the numbers of critical chains in the two cases are different, but the resulting M\"obius values become equal after cancellation.
\end{exa}

\subsection{Putting it all together}

We now have all the ingredients we need to proceed to the proof proper of Theorem~\ref{main:thm}

\bprf[Proof of Theorem~\ref{main:thm}]
As we saw in \eqref{sumoveremb:eq}, 
$
\mupstar(u,w) = \sum_\emb  \muembpstar{\emb}{w},
$
where the sum is over all embeddings $\emb$ of $u$ in $w$, so it suffices to show that 
\[
\muembpstar{\emb}{w} = \prod_{1 \leq j \leq \len{w}} \left\{ 
\begin{array}{ll} 
\mup(\emb(j), w(j)) +1 & \mbox{if $\emb(j)=0$ and $w(j-1)=w(j)$}, \\ 
\mup(\emb(j),w(j)) & \mbox{otherwise}.
\end{array} \right. 
\]
We will proceed by induction on $\len{w}=\ell$ with the result being trivially true when $\ell=1$.  

As we deduced from Lemma~\ref{descent:lem} and Corollary~\ref{lexdec:cor}, the critical chains $C$ from $w$ to $\emb$ must proceed by reducing $w(j)$ to $\emb(j)$ from right to left.  Let $k$ denote the largest index such that $w(k) \neq \emb(k)$.    A key consequence of 
Lemma~\ref{rightmost:lem} and
Proposition~\ref{msi:pr} is that the MSIs along $C$ as it reduces $w(k)$ to $\emb(k)$ depend on $w(k)$, $w(k-1)$ and $\emb(k)$, but not on the rest of $w$ and $\emb$.  In particular, the letters after position $k$ have no affect on the MSI structure, so it suffices to take $k=\ell$.   
Thus a critical chain $C$ from $w$ to $\emb$ must start by reducing $w(\ell)$ to $\emb(\ell)$, arriving at an embedding $\zeta$ with corresponding word $z$.  Next, $C$ must reduce $w(j)$ to $\emb(j)$ where $j < \ell$ is as large as possible with $w(j) \neq \emb(j)$.  We know that $\{z\}$ is then a 1-descent, and hence an MSI by Lemma~\ref{descent:lem}.  As in the proof of Lemma~\ref{cartesian:lem}, $C$ consists of any critical chain $C_1$ from $w$ to $\zeta$, followed by the MSI $\{z\}$, followed by any critical chain $C_2$ from the embedding $\zeta$ all the way to $\emb$.  Summing over such $C_1$ and $C_2$, Theorem~\ref{dmt:thm} gives
\bea
\muembpstar{\emb}{w} & = & \sum_{C_1, C_2} (-1)^{d(C)} \\
& = & \sum_{C_1, C_2} (-1)^{|\jint(C_1)|+|\jint(C_2)|} \\
& = & \left( \sum_{C_1} (-1)^{d(C_1)} \right) \left(\sum_{C_2} (-1)^{d(C_2)} \right)\\
& = & \muembpstar{\zeta}{w} \sum_{C_2} (-1)^{d(C_2)} .
\eea

The critical chains $C_2$ only reduce letters of $\zeta$ in positions $1, 2, \ldots, \ell-1$.  Again by Lemma~\ref{rightmost:lem} and Proposition~\ref{msi:pr}, the MSI structure will depend on the word $z' = z(1)\wdots z(\ell-1)$ and the embedding $\emb' = \emb(1)\wdots\emb(\ell-1)$ but not on $z(\ell)$ or $\emb(\ell)$.  Thus the critical chains $C_2$ and their $\jint$-intervals are in bijection with those of the interval $[u',z']$ that end at $\emb'$, where $u'$ is the word corresponding to $\emb'$.  Thus the sum over $C_2$ above equals $\muembpstar{\emb'}{z'}$ which, by induction on $\ell$, satisfies
\[
\muembpstar{\emb'}{z'} = \prod_{1 \leq j \leq \ell-1} \left\{ 
\begin{array}{ll} 
\mup(\emb(j), w(j)) +1 & \mbox{if $\emb(j)=0$ and $w(j-1)=w(j)$}, \\ 
\mup(\emb(j),w(j)) & \mbox{otherwise},
\end{array} \right. 
\]
since $z'(j) = w(j)$ for $1 \leq j \leq \ell-1$.  Since $\zeta(\ell) = \emb(\ell)$, it remains to show that 
\beq\label{muemb:eq}
\muembpstar{\zeta}{w} = \left\{ 
\begin{array}{ll} 
\mup(\zeta(\ell), w(\ell)) +1 & \mbox{if $\zeta(\ell)=0$ and $w(\ell-1)=w(\ell)$}, \\ 
\mup(\zeta(\ell),w(\ell)) & \mbox{otherwise}.
\end{array} \right.
\eeq

Since $\zeta$ and $w$ differ only in one position, we are in the situation of Proposition~\ref{msi:pr}.  So suppose first that $\zeta$ is rightmost and consider an interval $C(w',\zeta')$ on a maximal chain $C$ from $w$ to $\zeta$.  Then  $\zeta'$ is the rightmost embedding in $w'$ for the following reasons:  if $\zeta'(\ell)>_0 0$ then $\zeta'$ is the only embedding in $w'$; if $\zeta'(\ell)=0$ then $\zeta'=\zeta$, and $\zeta'$ is rightmost in $w'$ because it is in $w$ and $w'\le w$. So we can apply Proposition~\ref{msi:pr} with $w'$ and $\zeta'$ in place of $w$ and $\emb$.  By Condition~(a), we get that $C(w',\zeta')$ is an MSI if and only if $\czero(w'(\ell), \zeta'(\ell))$ is an MSI.  Therefore the MSIs on maximal chains $C$ from $w$ to $\zeta$ are in exact correspondence with MSIs on maximal chains $\czero$ from $w(\ell)$ to $\zeta(\ell)$.  Therefore, $C$ is a critical chain if and only if $\czero$ is a critical chain, and they are both covered by the same number of $\jint$-intervals.  We conclude that $\muembpstar{\zeta}{w} = \mup(\zeta(\ell), w(\ell))$ when $\zeta$ is rightmost, consistent with \eqref{muemb:eq}. 

Next suppose that $\zeta$ is not rightmost.  By Lemma~\ref{rightmost:lem}, we have $\zeta(\ell)=0$ and $w(\ell-1)\le_0 w(\ell)$.  Suppose first that $w(\ell-1)=w(\ell)$.
Consider an interval $C(w',\zeta')$ on a maximal chain $C$ from $w$ to $\zeta$.  If $\zeta'(\ell) >_0 0$, then $\zeta'$ is rightmost in $w'$, and we are in the case of the previous paragraph.  We get that $C(w',\zeta')$ is an MSI if and only if $\czero(w'(\ell), \zeta'(\ell))$ is an MSI.  This same equivalence applies if $w' < w$ in $\pstar$, since then $w'(\ell-1) =w(\ell-1)=w(\ell)>_0 w'(\ell)$ and so $\zeta'$ is rightmost in $w'$ by Lemma~\ref{rightmost:lem}.  So assume that $\zeta'=\zeta$ and $w'=w$.  Referring to Condition~(b$'$) of Proposition~\ref{msi:pr},  we see that since $w(\ell-1)=w(\ell)$, the condition that the open interval $\czero(w(\ell),\zeta(\ell))$ not contain an $x \geq_0 w(\ell-1)$ is automatically satisfied.   Hence $C(w,\zeta)$ is an MSI if and only if one of two mutually exclusive conditions is satisfied: $\czero(w(\ell),\zeta(\ell))$ is an MSI or $\czero(w(\ell),\zeta(\ell))$ is lexicographically first.  Putting this all together, we get that $C$ from $w$ to $\zeta$ is a critical chain if and only if $\czero$ is a critical chain or is lexicographically first.  In the former case, $C$ and $\czero$ are covered by the same number of $\jint$-intervals, while in the latter case $C(w,\zeta)$ is a single MSI.  This yields
\[
\muembpstar{\zeta}{w} = 
\mup(\zeta(\ell), w(\ell)) +1  \mbox{\ \ \ if $\zeta(\ell)=0$ and $w(\ell-1)=w(\ell)$}
\]
as in \eqref{muemb:eq}.

Finally, we suppose that $\zeta$ is not rightmost and that $w(\ell-1) <_0 w(\ell)$.  Consider an interval $C(w',\zeta')$ of a maximal chain $C$ from $w$ to $\zeta$.  Again Lemma~\ref{rightmost:lem} and Proposition~\ref{msi:pr} tell us that 
whether $C(w',\zeta')$ is an MSI does not depend on $w(1), \ldots, w(\ell-2)$ or $\zeta(1), \ldots, \zeta(\ell-2)$.  Therefore, it suffices to consider the case where $w$ has just two letters, i.e., $\ell=2$.  
We see that we are in exactly the situation of Corollary~\ref{special_case:co}, and we conclude that $\muembpstar{\zeta}{w} = \mup(0, w(\ell))$, as in \eqref{muemb:eq}.
\eprf

\begin{rem}
In the proof above, we could have proved the case when $\zeta$ is not rightmost and $w(\ell-1)=w(\ell)$ using the same technique used in the last paragraph for $w(\ell-1) <_0 w(\ell)$; 
Corollary~\ref{special_case:co} would give that $\muembpstar{\zeta}{w} = \mup(0, w(\ell)) +1$.  We chose to rely instead on Proposition~\ref{msi:pr} in order to make clear the connection between the lexicographically first chain and the $+1$ in Theorem~\ref{main:thm}.
\end{rem}

\section{Applications}
\label{a}

In this section we will show how the M\"obius function values for subword order, composition order and other special cases of generalized subword order mentioned in the Introduction all follow easily from Theorem~\ref{main:thm}.   First, however, we would like to prove a result about the homotopy type of certain $P^*$.

If $P$ is a finite poset and $x\in P$, then the {\it rank of $x$\/}, denoted $\rk(x)$, is the length of a
longest chain from a minimal element of $P$ to $x$.  In particular, minimal elements have rank 0.  The {\it rank of $P$\/} is
$$
\rk(P)=\max_{x\in P}\, \rk(x).
$$
For example, $P$ is an antichain if and only if $\rk(P)=0$.  
For $w \in P^*$, we will write $\rk(w)$ to denote the rank of $w$ in the interval $[\emptyset,w]$ of $P^*$.  Note that if $P$ is an antichain then $\rk(w)=|w|$ for $w \in P^*$.

Now
consider the {\it order complex\/}, $\Delta(x,y)$, of a finite interval $[x,y]$ in a poset $P$, which is the abstract simplicial complex consisting of all chains of $(x,y)$.  If $\Delta(x,y)$ has a topological property, we will also say that $[x,y]$ has the same property.  To prove Theorem~\ref{dmt:thm}, Babson and Hersh showed that $\Delta(x,y)$ is homotopic to a CW-complex with a cell for each critical chain and an extra  cell of dimension $0$.  The simplex in a critical chain $C$ giving rise to a critical cell is obtained by taking one element from each of the $\jint$-intervals and so has dimension $d(C)$.  This is all the information we need to prove the following result.
\bth
\label{homotopy}
Let $P$ be any finite poset with $\rk(P)\le 1$.  Then any interval $[u,w]$ in $P^*$ is homotopic to a wedge of $|\mu(u,w)|$ spheres all of dimension $\rk(w)-\rk(u)-2$.
\eth
\bprf
We claim that every MSI in 
a maximal chain of $[u,w]$
consists of one element.  Suppose, towards a contradiction, that there is a chain $C$ of the form~\ree{wordC} containing an MSI $I=\{v_i,v_{i+1},\ldots, v_k\}$ for some $i$ and $k>i$.  By Lemma~\ref{ascent:lem}, $I$ cannot contain an ascent.  And by Lemma~\ref{descent:lem}, $I$ cannot contain a 1-descent since otherwise it would contain a smaller SI.
So there must be some index $j$ so that only elements in position $j$ are decreased in passing from $v_{i-1}$ to $v_{k+1}$.  But then these elements form a chain of length at least 3 in $P_0$.  This contradicts the fact that the longest chain in $P$ has length at most 1.

Note that all the maximal chains in $[u,w]$ have length $\rk(w)-\rk(u)$.  And since any critical chain $C$ is covered by 1-element MSIs, the number of intervals in $\jint(C)$ is always $\rk(w)-\rk(u)-1$.  This also implies that there is no cancellation  in the sum of Theorem~\ref{dmt:thm}.  So $\mu(u,w)$ is, up to sign, the number of critical chains.
  It follows that the CW-complex discussed above must be constructed from a 0-cell together with $|\mu(u,w)|$ cells of dimension $\rk(w)-\rk(u)-2$.  The only way to construct such a complex is as given in the statement of the theorem.
\eprf

We note that Babson and Hersh~\cite{bh:dmf} show that if every MSI of an interval $[x,y]$ in a poset   is a singleton, then $[x,y]$ is shellable and hence a wedge of spheres.  (They assume that the interval is pure, but their proof goes through for non-pure posets.)  So our proof above actually shows that $[u,w]$ is shellable whenever $\rk(P)\le 1$.

\subsection{Bj\"orner's formula for subword order}

Subword order on the alphabet $A = \{1, 2, \ldots, s\}$ corresponds to the case when $\poset$ is an antichain with elements $A$.   Let us determine what Theorem~\ref{main:thm} yields in this case.  For $u, w \in A^*$, suppose $\emb$ is an embedding of $u$ in $w$.  For each $j$ with $1 \leq j \leq \len{w}$, there are two cases.  The first is that $\emb(j) = w(j) \neq 0$, in which case $\mup(\emb(j),w(j))=1$.  The more interesting situation is when $\emb(j)=0$, in which case $\mup(\emb(j), w(j))=-1$.  So, if $w(j-1)=w(j)$, then $\emb$ contributes 0 to the sum in Theorem~\ref{main:thm}.  Thus we can restrict to \emph{normal embeddings}, meaning that $\emb(j)\neq 0$ whenever $w(j-1)=w(j)$.  For example, if $w=1122121$ and $u=121$, then there are exactly two normal embeddings of $u$ in $w$, namely $0102100$ and $0102001$.  Let us denote the number of normal embeddings of $u$ in $w$ by $\binom{w}{u}_n$.  Putting these observations together, we get Bj\"orner's result from Theorem~\ref{main:thm}.

\bth[\cite{bjo:mfs}]  If $u, w \in A^*$, then 
\[
\mu(u,w) = (-1)^{\len{w}-\len{u}} \binom{w}{u}_n.  
\] 
\eth

In the same paper, Bj\"orner also derived the homotopy type of $[u,w]$, and this result follows immediately from Theorem~\ref{homotopy}.
\bth[\cite{bjo:mfs}]  If $u, w \in A^*$, then $[u,w]$ is homotopic to a wedge of $\binom{w}{u}_n$ spheres, 
all of dimension $|w|-|u|-2$.\hqed
\eth

\subsection{Generalized subword order for rooted forests}

We now consider the generalization of Bj\"orner's result to rooted forests given by Sagan and Vatter~\cite{sv:mfc}.
Clearly,  $P$ is a rooted forest if and only if every element $x\in P_0-\{0\}$ covers exactly one element, denoted $x^-$, of $\pzero$.  We will show how Theorem~\ref{main:thm} gives the formula for $\mupstar$ as stated in~\cite{sv:mfc}; their statement for the M\"obius function of composition order is almost identical and follows immediately.

For $\poset$ a rooted forest and $u, w \in P^*$, let $\emb$ be an embedding of $u$ in $w$.  Note that for $x, y \in \pzero$,
\[
\mup(x,y)=\left\{ \begin{array}{cl}
+1 & \mbox{if $y=x$},\\
 -1 & \mbox{if $y \coversinp x$},\\
 0 & \mbox{otherwise.}
 \end{array}
 \right.
\]
Therefore, if $\emb(j)\neq 0$, for $\emb$ to contribute a nonzero amount to the sum in Theorem~\ref{main:thm}, there are two possibilities:
\begin{itemize}
\item $w(j)=\emb(j)$, which will contribute a $1$ to the product, or
\item  $w(j) \coversinp \emb(j)$, which will contribute a $-1$.  
\end{itemize}
If $\emb(j)=0$, there are also two possibilities that will allow $\emb$ to have a nonzero contribution: 
\begin{itemize}
\item $w(j)$ is a minimal element of $P$ and $w(j-1)\neq w(j)$, which will contribute a $-1$ to the product, or
\item $w(j)$ is not minimal and $w(j-1)=w(j)$, which will contribute a $1$. 
\end{itemize} 
These four conditions on $w$ and $\emb$ can be seen to be equivalent to those on the generalized version of \emph{normal embedding} in \cite{sv:mfc}, defined there as an embedding $\emb$ of $u$ in $w$ satisfying the following two conditions.
\begin{enumerate}
\item For $1 \leq j \leq \len{w}$ we have $\emb(j) = w(j), w(j)^-$ or $0$.  
\item For all $x \in \poset$ and every run $[r,t]$ of $x$'s in $w$, we have
\begin{enumerate}
\item $\emb(j) \neq 0$ for all $j$ with $r < j \leq t$ if $x$ is minimal in $P$, 
\item $\emb(r) \neq 0$ otherwise.
\end{enumerate}
\end{enumerate}
The \emph{defect} $\defect(\emb)$ of a normal embedding $\emb$ of $u$ in $w$ is defined in \cite{sv:mfc} as
\[
\defect(\emb) = \#\{i : \emb(i) = w(i)^-\}.
\]
Referring to Theorem~\ref{main:thm}, we see that the defect is exactly the number of $j$'s, $1\leq j \leq \len{w}$, that will contribute $-1$ to a nonzero product, while all other $j$'s in a normal embedding will contribute $+1$.  Putting this all together, we get \cite[Theorem~6.1]{sv:mfc}.

\bth[\cite{sv:mfc}]
Let $\poset$ be a rooted forest.  Then the M\"obius function of $\pstar$ is given by
\[
\mupstar(u, w) = \sum_\emb (-1)^{\defect(\emb)},
\]
where the sum is over all normal embeddings $\emb$ of $u$ in $w$.  
\eth

Restricting to the composition poset, which arises when $P=\bbP$, everything stays the same, except that we can write $w(j)-1$ in place of $w(j)^-$ and we can change the language in the second condition on a normal embedding to read:
\begin{enumerate}
\item[2.]  For all $k \geq 1$ and every run $[r,t]$ of $k$'s in $w$, we have
\begin{enumerate}
\item $\emb(j) \neq 0$ for all $j$ with $r < j \leq t$ if $k=1$, 
\item $\emb(r) \neq 0$ if $k \geq 2$.
\end{enumerate}
\end{enumerate}

\subsection{Connection with Chebyshev polynomials}\label{chebyshev:ss}
As promised, a connection between generalized subword order and Chebyshev polynomials follows easily from Theorem~\ref{main:thm}.  Consider the poset $\Lambda$ from Figure~\ref{Lambda:fig} and the intervals $[1^i, 3^j]$ in $\Lambda^*$.  To describe the corresponding M\"obius function values, consider the Chebyshev polynomials $T_n(x)$ of the first kind, which can be defined recursively by $T_0(x)=1$, $T_1(x)=x$, and 
\beq\label{chebrecursive:eq}
T_n(x) = 2x\, T_{n-1}(x) - T_{n-2}(x)
\eeq
for $n>1$.  An equivalent definition which is more suitable for our purposes is obtained by replacing \eqref{chebrecursive:eq} by
\beq\label{chebyshev:eq}
T_n(x) = \frac{n}{2} \sum_{k=0}^{\lfloor \frac{n}{2}\rfloor} \frac{(-1)^k}{n-k}\binom{n-k}{k} (2x)^{n-2k}
\eeq
for $n>1$.  One consequence of either definition is that the coefficient of $x^m$ in $T_n(x)$, which we will denote by $\spn{x^m}T_n$, is nonzero only if $m$ and $n$ have the same parity.  The following result, which was conjectured in \cite{sv:mfc} and first proved in \cite{tom:gcp}, concerns such coefficients. 

\bth[\cite{tom:gcp}] \label{chebyshev:thm} Considering intervals in $\Lambda^*$,
for all $0 \leq i \leq j$, 
$$
\mulambdastar(1^i, 3^j)=\spn{x^{j-i}}T_{i+j}(x).
$$  
\eth

\bprf
First, we check the result for $i+j=0$.
We must have $i=j=0$ and $[1^i,3^j]$ is a single element poset, consistent with $T_0(x)=1$.

Otherwise, $j\ge1$ and there are $\binom{j}{i}$ embeddings $\emb$ of $1^i$ in $3^j$, of which $\binom{j-1}{i}$ satisfy $\emb(1)=0$ while $\binom{j-1}{i-1}$ satisfy $\emb(1)=1$ (where binomial coefficients of the form $\binom{n}{k}$ with $k<0$ are considered 0 as usual).  By Theorem~\ref{main:thm}, the former type of embeddings each contribute $(-1)^i2^{j-i-1}$ to the M\"obius function, while the latter type each contribute $(-1)^i 2^{j-i}$.  Thus 
\bea
\mulambdastar(1^i,3^j) & = & (-1)^i 2^{j-i-1} \left( \binom{j-1}{i} + 2\binom{j-1}{i-1}\right) \\
 & = & (-1)^i 2^{j-i-1} \left(\binom{j-1}{i-1} + \binom{j}{i} \right) \\
 & = & (-1)^i 2^{j-i-1} \frac{i+j}{j} \binom{j}{i}.
\eea
This last expression is now readily checked to be the coefficient of $x^{j-i}$ when $n=i+j$ in \eqref{chebyshev:eq}.
\eprf

Although Tomie did not derive the homotopy type for these intervals in $\La^*$, we obtain the information easily from Theorem~\ref{homotopy}.
\bth 
For all $0 \leq i \leq j$, the interval $[1^i, 3^j] $ in $\Lambda^*$
is homotopic to a wedge of $|\spn{x^{j-i}}T_{i+j}(x)|$ spheres, 
all of dimension $2j-i-2$.\hqed
\eth

\subsection{Tomie's generalized Chebyshev polynomials}

The main result of \cite{tom:gcp} is more general than Theorem~\ref{chebyshev:thm}.  For $s\geq 1$, Tomie considers the poset, which we denote by $\ptomie$, that consists of an $s$-element antichain $\{1, 2, \dots, s\}$ with a top element $s+1$ added.  Letting $s=2$ gives the poset $\Lambda$.  Along the same lines, Tomie recursively defines generalized Chebyshev polynomials $T^s_n(x)$, and gives a closed-form expression for the coefficients of $T^s_n(x)$ which, after a change of variables, can be written as
\beq\label{tomie:eq}
T^s_n(x) = \sum_{k=0}^{\lfloor\frac{n}{2}\rfloor} (-1)^k s^{n-2k-1}\left( \binom{n-k}{k} s - \binom{n-k-1}{k} \right) x^{n-2k}.
\eeq
for $n\geq 0$ and $s\geq 1$. 

The main result of \cite{tom:gcp} again follows as a special case of Theorem~\ref{main:thm}, as we now show.

\bth[\cite{tom:gcp}]
Considering intervals in $(\ptomie)^*$,
for all $0 \leq i \leq j$ and $s \geq 1$,
$$
\mu(1^i, (s+1)^j)=\spn{x^{j-i}}T^s_{i+j}(x).
$$
\eth

\bprf
From \eqref{tomie:eq}, we get that
\[
T^s_{i+j}(x) = \sum_{k=0}^{\lfloor\frac{i+j}{2}\rfloor} (-1)^k s^{i+j-2k-1} \left( \binom{i+j-k}{k} s - \binom{i+j-k-1}{k} \right) x^{i+j-2k}.
\]
Considering the term in the sum where $k=i$, we get that
the coefficient of $x^{j-i}$ in $T^s_{i+j}(x)$ is 
\[
(-1)^i s^{j-i-1} \left( \binom{j}{i} s - \binom{j-1}{i} \right),
\]
which equals
\beq\label{gencheb:eq}
(-1)^i s^{j-i-1} \left( \binom{j-1}{i-1} s +  \binom{j-1}{i}(s-1) \right)
\eeq
whenever $j \geq 1$.

Now consider $\mu(1^i, (s+1)^j)$ as determined by Theorem~\ref{main:thm}.  
When $i+j=0$, we must have $i=j=0$ and $[1^i,(s+1)^j]$ is a single element poset, consistent with $T^s_0(x)=1$ from \eqref{tomie:eq}.  Otherwise, $j\ge1$ and there are $\binom{j}{i}$ embeddings $\emb$ of $1^i$ in $(s+1)^j$, of which $\binom{j-1}{i}$ satisfy $\emb(1)=0$ while $\binom{j-1}{i-1}$ satisfy $\emb(1)=1$.  By Theorem~\ref{main:thm}, the former type of embeddings each contribute $(-1)^i (s-1) s^{j-i-1}$ to the M\"obius function, while the latter type each contribute $(-1)^i s^{j-i}$.  Thus $\mu(1^i, (s+1)^j)$ equals the expression~\eqref{gencheb:eq}, as required.
\eprf

Using Theorem~\ref{homotopy} one last time, we obtain the following result.
\bth
For all $0 \leq i \leq j$, the interval $[1^i, (s+1)^j] $ in $(\ptomie)^*$
is homotopic to a wedge of $|\spn{x^{j-i}}T_{i+j}^s(x)|$ spheres, 
all of dimension $2j-i-2$.\hqed
\eth

\section{Closing Remarks}
\label{cr}

There has also been interest in \emph{generalized factor order} on $P^*$ which is defined like generalized subword order except that one requires the indices $i_1,i_2,\ldots,i_k$ to be consecutive in~\ree{gso:def}.  Bj\"orner \cite{bjo:mff} found a recursive formula for the M\"obius function in the case of \emph{ordinary factor order}, which corresponds to $P$ being an antichain.  In particular, he showed that the only possible M\"obius values are $0,\pm1$ and that the order complex of every interval is homotopic to either a ball or a sphere.  In his thesis, see \cite{wil:mfg}, Willenbring reproved Bj\"orner's results in an elucidating way using critical chains and found a more general formula for rooted forests.  The latter is quite complicated.  It would be very interesting if one could find a simpler formula more along the lines of Theorem~\ref{main:thm}.

The analogue of generalized factor order for $\fS$ is called the consecutive pattern poset.    Somewhat surprisingly (given the fact that Wilf's question for ordinary patterns has not been fully answered), Bernini, Ferrari, and Steingr\'{\i}msson \cite{bfs:mfc} gave a complete description of the M\"obius function in the consecutive case.  Even more surprisingly, Sagan and Willenbring \cite{sw:cpp} were able to give a proof of this result using critical chains which closely parallels the one Willenbring gave for factor order of an antichain.  This led them to define a sequence of partial orders on $\bbP^*$, denoted $P_k$ for $k=0,1,2,\dots,\infty$, where $P_0$ is ordinary factor order, $P_\infty$ contains consecutive pattern order as a convex subposet, and every $P_k$ has essentially that same M\"obius function.  So this sequence of interpolating posets gives an explanation of the coincidence noted above.

\

\noindent
{\it Acknowledgements.}  
The authors thank Volker Strehl and the anonymous referees for helpful comments.
Calculations were performed using John Stembridge's posets package \cite{ste:spc}.

\bibliographystyle{plain}
\bibliography{gso}

\begin{thebibliography}{10}

\bibitem{bh:dmf}
Eric Babson and Patricia Hersh.
\newblock Discrete {M}orse functions from lexicographic orders.
\newblock {\em Trans. Amer. Math. Soc.}, 357(2):509--534 (electronic), 2005.

\bibitem{bfs:mfc}
Antonio Bernini, Luca Ferrari, and Einar Steingr{\'{\i}}msson.
\newblock The {M}\"obius function of the consecutive pattern poset.
\newblock {\em Electron. J. Combin.}, 18(1):Paper 146, 2011.

\bibitem{bjo:mfs}
Anders Bj{\"o}rner.
\newblock The {M}\"obius function of subword order.
\newblock In {\em Invariant theory and tableaux ({M}inneapolis, {MN}, 1988)},
  volume~19 of {\em IMA Vol. Math. Appl.}, pages 118--124. Springer, New York,
  1990.

\bibitem{bjo:mff}
Anders Bj{\"o}rner.
\newblock The {M}\"obius function of factor order.
\newblock {\em Theoret. Comput. Sci.}, 117(1-2):91--98, 1993.

\bibitem{bjjs:mfs}
Alexander Burstein, V\'{\i}t Jel\'{\i}nek, Eva Jel\'{\i}nkov\'a, and Einar
  Steingr\'{\i}msson.
\newblock The {M}\"obius function of separable and decomposable permutations.
\newblock {\em J. Combin. Theory Ser. A}, 118(8):2346--2364, 2011.

\bibitem{for:dmt}
Robin Forman.
\newblock A discrete {M}orse theory for cell complexes.
\newblock In {\em Geometry, topology, \& physics}, Conf. Proc. Lecture Notes
  Geom. Topology, IV, pages 112--125. Internat. Press, Cambridge, MA, 1995.

\bibitem{for:mtc}
Robin Forman.
\newblock Morse theory for cell complexes.
\newblock {\em Adv. Math.}, 134(1):90--145, 1998.

\bibitem{for:ugd}
Robin Forman.
\newblock A user's guide to discrete {M}orse theory.
\newblock {\em S\'em. Lothar. Combin.}, 48:Art.\ B48c, 35, 2002.

\bibitem{kru:twq}
Joseph~B. Kruskal.
\newblock The theory of well-quasi-ordering: {A} frequently discovered concept.
\newblock {\em J. Combinatorial Theory Ser. A}, 13:297--305, 1972.

\bibitem{sv:mfc}
Bruce~E. Sagan and Vincent Vatter.
\newblock The {M}\"obius function of a composition poset.
\newblock {\em J. Algebraic Combin.}, 24(2):117--136, 2006.

\bibitem{sw:cpp}
Bruce~E. Sagan and Robert Willenbring.
\newblock Discrete {M}orse theory and the consecutive pattern poset.
\newblock {\em J. Algebraic Combin.}
\newblock To appear.
  \href{http://arxiv.org/abs/1107.3262}{\texttt{arxiv:1107.3262}}.

\bibitem{sta:ec1}
Richard~P. Stanley.
\newblock {\em Enumerative combinatorics. {V}ol. 1}, volume~49 of {\em
  Cambridge Studies in Advanced Mathematics}.
\newblock Cambridge University Press, Cambridge, 1997.
\newblock With a foreword by Gian-Carlo Rota. Corrected reprint of the 1986
  original.

\bibitem{st:mfp}
Einar Steingr{\'{\i}}msson and Bridget~Eileen Tenner.
\newblock The {M}\"obius function of the permutation pattern poset.
\newblock {\em J. Comb.}, 1(1):39--52, 2010.

\bibitem{ste:spc}
John~R. Stembridge.
\newblock {SF}, posets and coxeter/weyl.
\newblock \newline Available from
  \href{http://www.math.lsa.umich.edu/~jrs/maple.html}{\url{www.math.lsa.umich%
.edu/~jrs/maple.html}}.

\bibitem{tom:gcp}
Masaya Tomie.
\newblock A generalization of the {C}hebyshev polynomials and nonrooted posets.
\newblock {\em Int. Math. Res. Not. IMRN}, (5):856--881, 2010.

\bibitem{wil:pp}
Herbert~S. Wilf.
\newblock The patterns of permutations.
\newblock {\em Discrete Math.}, 257(2-3):575--583, 2002.

\bibitem{wil:mfg}
Robert Willenbring.
\newblock The {M}\"obius function of generalized factor order.
\newblock Preprint.

\end{thebibliography}
\end{document}